\documentclass[a4paper, 12pt]{article}
\usepackage[OT2, T1]{fontenc}
\usepackage{cite}
\usepackage{amsmath}
\usepackage{amsfonts}
\usepackage{amssymb}
\usepackage{amsxtra}
\usepackage{mathrsfs}
\usepackage{mathtools}
\usepackage{tensor}
\usepackage{tikz-cd}
\usepackage[all]{xy}
\usepackage[thmmarks, amsmath, amsthm]{ntheorem}
\usepackage[pagebackref=true]{hyperref}
\hypersetup{hidelinks}

\usepackage{graphicx}
\usepackage[new]{old-arrows} 

\usepackage[portrait, top=2cm, bottom=2cm, left=2cm, right=2cm]{geometry}
\usepackage{enumerate}

\numberwithin{equation}{section}
\renewcommand{\labelenumi}{\textup{(\roman{enumi})}}
\renewcommand{\theenumi}{\textup{(\roman{enumi})}}

\theoremstyle{plain}
\newtheorem{thm}{Theorem}[section]
\newtheorem{lemm}[thm]{Lemma}
\newtheorem{prop}[thm]{Proposition}
\newtheorem{coro}[thm]{Corollary}

\theoremstyle{definition}

\newtheorem{exam}[thm]{Example}

\newcommand{\ab}{\operatorname{ab}}
\newcommand{\Aut}{\operatorname{Aut}}
\newcommand{\Br}{\operatorname{Br}}
\newcommand{\CH}{\operatorname{CH}}
\newcommand{\Cor}{\operatorname{Cor}}
\newcommand{\Der}{\operatorname{Der}}
\newcommand{\End}{\operatorname{End}}
\newcommand{\et}{\mathrm{\acute{e}t}}
\newcommand{\ev}{\operatorname{ev}}
\newcommand{\Ext}{\operatorname{Ext}}
\newcommand{\Gal}{\operatorname{Gal}}
\newcommand{\GL}{\mathbf{GL}}
\newcommand{\gp}{\textup{-gp}}
\newcommand{\Hom}{\operatorname{Hom}}
\newcommand{\id}{\operatorname{id}}
\newcommand{\Img}{\operatorname{Im}}
\newcommand{\Inn}{\operatorname{Inn}}
\newcommand{\Inv}{\operatorname{Inv}}
\newcommand{\inv}{\operatorname{inv}}
\newcommand{\Ker}{\operatorname{Ker}}
\newcommand{\Mat}{\mathrm{Mat}}

\newcommand{\nr}{\operatorname{nr}}
\newcommand{\Pic}{\operatorname{Pic}}
\newcommand{\Rad}{\operatorname{Rad}}
\newcommand{\rank}{\operatorname{rank}}
\newcommand{\SL}{\mathbf{SL}}
\newcommand{\Spec}{\operatorname{Spec}}

\newcommand{\tg}{\operatorname{tg}}
\newcommand{\tors}{\operatorname{tors}}

\newcommand{\Xbb}{\mathbb{X}}

\newcommand{\Cbf}{\mathbf{C}}
\newcommand{\Fbf}{\mathbf{F}}
\newcommand{\Gbf}{\mathbf{G}}

\newcommand{\Pbf}{\mathbf{P}}
\newcommand{\Qbf}{\mathbf{Q}}
\newcommand{\Rbf}{\mathbf{R}}
\newcommand{\Zbf}{\mathbf{Z}}

\newcommand{\Acal}{\mathcal{A}}

\newcommand{\Fcal}{\mathcal{F}}
\newcommand{\Ocal}{\mathcal{O}}
\newcommand{\Rcal}{\mathcal{R}}

\renewcommand{\H}{\mathrm{H}}
\newcommand{\M}{\mathrm{M}}
\newcommand{\N}{\mathrm{N}}
\newcommand{\Z}{\mathrm{Z}}

\title{Zero-cycles and the algebraic Brauer--Manin obstruction on homogeneous spaces}
\author{ Nguy\~{\^{e}}n M\d{a}nh Linh \\ \vspace*{-1ex} \small \textit{Universit{\'e} d'Artois} \and Tianyang Wang \\ \vspace*{-1ex} \small \textit{Sorbonne Universit{\'e}}}
\date{\today}

\begin{document}
\maketitle

\begin{abstract}
    For smooth proper models of homogeneous spaces of connected linear algebraic groups over number fields, we prove that the algebraic Brauer group controls the existence of zero-cycles of prescribed degree. This refines a theorem of Harpaz and Wittenberg, namely Conjecture ($E_1$) for such varieties. The result contrasts with the situation for rational points and with weak approximation for both rational points and zero-cycles, where the transcendental Brauer--Manin obstruction is strictly indispensable. Our proof builds on the theorem of Harpaz and Wittenberg on zero-cycles on homogeneous spaces and combines a reduction lemma for nonabelian Galois cohomology due to Chernousov--Gille--Reichstein with Demu\v{s}kin's description of the maximal pro-$p$ Galois groups of local fields.
\end{abstract}


{\em 2020 Mathematics subject classifications (MSC).} 11R34, 14G12, 14M17.

\tableofcontents

\section{Introduction} \label{sec:Intro}

Let $X$ be a smooth proper geometrically integral variety over a number field $k$. Manin's foundational idea \cite{Manin1971Brauer} associates with $X$ the {\em Brauer--Manin pairing}
\[
    \prod_{v\in\Omega_k} X(k_v) \times \Br(X) \longrightarrow \Qbf/\Zbf, \qquad ((P_v)_v,\Acal) \longmapsto \sum_{v\in\Omega_k} \inv_v \bigl(\Acal(P_v)\bigr),
\]
where $\Omega_k$ denotes the set of places of $k$, where $\Br(X)$ is the Brauer group of $X$, and where
\[
    \inv_v = \inv_{k_v} \colon \Br(k_v) \longhookrightarrow \Qbf/\Zbf
\]
are the invariant maps from local class field theory. Without the properness assumption on $X$, one replaces $\Br(X)$ by the {\em unramified} Brauer group $\Br_{\nr}(X)$. The global reciprocity law (Albert--Brauer--Hasse--Noether) shows that the diagonal image of $X(k)$ is orthogonal to $\Br_{\nr}(X)$, and that every family $(P_v)_v$ of local points is orthogonal to the subgroup
\[
    \Br_0(X) \coloneq \Img \bigl(\Br(k) \longrightarrow \Br(X)\bigr)
\]
of {\em constant} Brauer classes. In particular, the vacuity of the {\em Brauer--Manin set} 
\[
    \left\{(P_v)_v \in \prod_{v \in \Omega_k} X(k_v) \mid \sum_{v \in \Omega_k} \inv_v\bigl(\Acal(P_v)\bigr) = 0 \text{ for all } \Acal \in \Br_{\nr}(X)\right\}
\]
is an obstruction to the Hasse principle for $X$. By continuity, the Brauer--Manin pairing also gives an obstruction to weak approximation.

A conjecture of Colliot-Th\'el\`ene predicts that the Brauer--Manin obstruction is the only obstruction to the Hasse principle and weak approximation on smooth {\em rationally connected} varieties over number fields; see \cite[Conjecture~14.1.2]{CTS2021Brauer}. Homogeneous spaces of connected linear algebraic groups (or, more generally, geometrically unirational varieties) belong to this class. For principal homogeneous spaces under connected linear algebraic groups without factors of type $E_8$, an affirmative answer follows from Sansuc's work \cite[Corollaires~8.7 and~8.13]{Sansuc1981Arithmetique}. Borovoi proved it for homogeneous spaces with connected geometric stabilizers and, when the ambient group is simply connected, with commutative stabilizers \cite[Theorems~2.2 and~2.3]{Borovoi1996BrauerManin}. For homogeneous spaces of semisimple simply connected groups, Harpaz and Wittenberg established the conjecture for finite supersolvable geometric stabilizers in the sense of \cite[D\'efinition~6.4 and Th\'eor\`eme~B]{HW2020Galois}, while Demeio recently settled the case of finite solvable geometric stabilizers \cite[Theorem~1.1]{Demeio2026Solvable}. The case of arbitrary {\em finite noncommutative} stabilizers nevertheless remains open, and it is the essential remaining case according to the reduction theorem of Demarche and Lucchini Arteche \cite[Th\'eor\`eme~1.1]{DLA2019Reduction}.

The distinction between the full unramified Brauer group and its algebraic part
\[
    \Br_{\nr,1}(X)\coloneq
    \Ker\bigl(\Br_{\nr}(X)\longrightarrow\Br_{\nr}(X_{\bar{k}})\bigr)
\]
is crucial. For the homogeneous spaces in Borovoi's theorem, the algebraic Brauer group already controls the Hasse principle and weak approximation. This is no longer true for finite noncommutative stabilizers. Demarche, Lucchini Arteche, and Neftin produced the first examples of {\em transcendental} Brauer--Manin obstructions to weak approximation on homogeneous spaces of $\SL_n$ with finite constant $p$-primary stabilizers \cite[Corollary~5.3 and Example~5.4]{DLAN2017Grundwald}. Further examples and explicit computations were obtained independently by Demeio \cite[Theorem~1.5]{Demeio2026Descent} and Lagarde \cite[Theorem~1.5]{Lagarde2025Unramified}. These works show that $\Br_{\nr,1}(X)$ does not in general account for all Brauer--Manin obstructions to weak approximation. The first-named author recently constructed homogeneous spaces $X$ for which $\Br_{\nr,1}(X)=\Br_0(X)$ and which are counterexamples to the Hasse principle \cite[Theorem~1]{Linh2026Insufficiency}. This confirms that transcendental Brauer classes can also be indispensable for the existence of rational points.

The Brauer--Manin formalism extends from rational points to {\em zero-cycles} ({\em i.e.}, finite formal sums of closed points). If $V$ is a smooth geometrically integral $k$-variety, evaluation at closed points followed by corestriction gives rise to local pairings
\[
    \Z_0(V_{k_v}) \times \Br_{\nr}(V) \longrightarrow \Qbf/\Zbf \qquad (v \in \Omega_k)
\]
which assemble into a global pairing. If $V$ is proper, these pairings factor through the local Chow groups $\CH_0(V_{k_v})$ of rational equivalence. {\em Conjecture \textup{(E)}} of Colliot-Th\'el\`ene, Sansuc, Kato, and Saito \cite{CTS1981Chow,KatoSaito1986Global,CT1995Arithmetic}, in the formulation of van Hamel and Wittenberg \cite{vanHamel2003Brauer,Wittenberg2012ZeroCycles}, predicts the sufficiency of the Brauer--Manin obstruction to the existence of zero-cycles of degree one and to Chow-group weak approximation for zero-cycles on smooth proper geometrically integral varieties. According to \cite[Rappel~1.3]{HW2020Galois}, for (smooth proper) rationally connected varieties $V/k$, it asserts the exactness of
    \begin{equation} \label{eq:ExactSequenceConjectureE} \tag{E}
        \CH_0(V) \longrightarrow \prod_{\substack{ v \in \Omega_k \\ v \mid \infty}} \frac{\CH_0(V_{k_v})}{\N_{\bar{k}_v/k_v}\bigl(\CH_0(V_{\bar{k}_v})\bigr)} \times \prod_{\substack{ v \in \Omega_k \\ v \nmid \infty}} \CH_0(V_{k_v}) \longrightarrow \Hom\bigl(\Br(V),\Qbf/\Zbf\bigr).
    \end{equation}
Harpaz and Wittenberg settled Conjecture \textup{(E)} for smooth compactifications of homogeneous spaces of connected linear algebraic groups \cite[Th\'eor\`eme~A]{HW2020Galois}. This was later generalized by the first-named author to a descent theorem for zero-cycles \cite[Theorem~4.1]{Linh2026Descent}.

In the above context of zero-cycles, Harari raised the question of whether we can replace the full Brauer group $\Br(V)$ in exact sequence \eqref{eq:ExactSequenceConjectureE} with its algebraic subgroup $\Br_1(V)$. The preceding examples for rational points of Demarche--Lucchini Arteche--Neftin, Demeio, and Lagarde already give a negative answer. Indeed, the obstructed families of local points define families of local zero-cycles of degree one which are orthogonal to the algebraic Brauer group but not to the full unramified Brauer group. Since sequence \eqref{eq:ExactSequenceConjectureE} is a complex, such families cannot be approximated by global zero-cycles. The situation is strikingly different if one only asks for the existence of a zero-cycle of a prescribed degree. The aim of the present article is to prove that for homogeneous spaces of connected linear algebraic groups, the algebraic unramified Brauer group gives the only obstruction to this existence problem. Here is our main result.

\begin{thm} \label{thm:Main}
    Let $X$ be a homogeneous space of a connected linear algebraic group over a number field $k$, and let $\delta \in \Zbf$. If there exists a family of degree $\delta$ local zero-cycles
    \[
        z_\Omega \in \prod_{v \in \Omega_k} \Z_0(X_{k_v})
    \]
    that is orthogonal to $\Br_{\nr,1}(X)$, then there exists a zero-cycle of degree $\delta$ on $X$.
\end{thm}

Note that the above ``equidegree'' requirement is not restrictive. Indeed, according to Lemma~\ref{lemm:OrthogonalityToBr0} below, if $z_\Omega = (z_v)_v$ is any family orthogonal to $\Br_0(X)$, then the coordinates $z_v$, at the finite places, have a common degree $\delta \in \Zbf$. Furthermore, we can modify the coordinates $z_v$ at the archimedean places by adding suitable multiples of zero-cycles of the form $\N_{\bar{k}_v/k_v}(\bar{P}_v)$, for $\bar{P}_v \in X_{\bar{k}_v}(\bar{k}_v)$ (these contribute trivially to the Brauer--Manin pairing, by functoriality and the vanishing of $\Br(\bar{k}_v)$), in such a way that $\deg(z_v) = \delta$ for all $v$.

This refinement avoids the need to determine the full unramified Brauer group, which is often difficult to compute. The algebraic part can instead be studied through the identification $\Br_{\nr,1}(X)/\Br_0(X)\simeq\H^1(k,\Pic(X^c_{\bar{k}}))$, where $X^c$ is a smooth compactification of $X$.

\medskip

The manuscript is organized as follows. 

Section~\ref{sec:Preliminaries} fixes the conventions on zero-cycles, Brauer groups, and the Brauer--Manin pairing that are used throughout the proof. We also record a result of Chernousov, Gille, and Reichstein, which enables a reduction of the calculation to homogeneous spaces with finite stabilizers. When the acting group has trivial geometric Picard group, Proposition~\ref{prop:BrNrSubgroupOfH2} identifies the unramified geometric Brauer group of such a homogeneous space with a subgroup of the second cohomology of its finite geometric stabilizer with coefficients in $\Qbf/\Zbf$. The latter is canonically isomorphic to the Pontryagin dual of the Schur multiplier, allowing explicit computations.

In Section~\ref{sec:ExistenceOfZeroCycles}, we proceed to prove Theorem~\ref{thm:Main}. Our idea is as follows. Assume the existence of a family $(z_v)_v$ of local zero-cycles
orthogonal to $\Br_{\nr,1}(X)$. We ``correct'' the $p$-primary part of the transcendental Brauer--Manin defect on this family by modifying one coordinate $z_v$ for some $v$ lying above $p$. To this end, we first apply the reduction from Section~\ref{sec:Preliminaries} to pass to finite stabilizers, and then pass to a Sylow $p$-subgroup $S$. The key ingredient in the proof is a local flexibility statement over a sufficiently large $p$-adic extension $K$, namely, that any prescribed functional $\theta$ on $\H^2(S,\Qbf/\Zbf)$ can be realized via pullback along a homomorphism $\rho\colon \Gal(\bar{K}/K) \to S$. After replacing the acting group by a suitable cover, we can turn $\rho$ into a local point on the homogeneous space, whose evaluation is precisely $\theta$. We then obtain the desired correcting zero-cycle by pushing forward the difference between this point and a base point. This local flexibility lemma relies on Demu\v{s}kin's description of local $p$-maximal Galois groups.

In Section~\ref{sec:Example}, we illustrate some aspects of the main theorem through various numerical examples. We first give an example for which the algebraic Brauer--Manin obstruction in Theorem~\ref{thm:Main} is optimal, in the sense that $\Br_{\nr,1}(X)/\Br_0(X)$ cannot be replaced by a proper subgroup. We next generalize the construction in \cite{Linh2026Insufficiency} to obtain new counterexamples to the Hasse principle on homogeneous spaces with finite {\em perfect} geometric stabilizers (so the algebraic Brauer group consists of constant classes). Our refinement is particularly useful here: since the algebraic Brauer group consists of constant classes, the existence of points over every completion of $k$ already guarantees the existence of a global zero-cycle of degree one. We then construct explicitly these zero-cycles.
\section{Preliminaries} \label{sec:Preliminaries}

\subsection{Notations and conventions} \label{subsec:Notations}

The following notations and conventions shall be used throughout the article.

\par\medskip
\noindent{\bf Fields}\par

Unless otherwise stated, all fields are of characteristic zero. For an algebraic object $X$ (such as a variety, a zero-cycle, or a cohomology class) defined over $K$, we denote its base change to a field extension $L/K$ by $X_L$ or $X|_L$. An algebraic closure of $K$ is denoted by $\bar{K}$. 

If $k$ is a number field, we denote by $\Omega_k$ the set of its places. For $v \in \Omega_k$, denote by $k_v$ the corresponding completion, and by $\inv_v$ or $\inv_{k_v}$ the local invariant map 
\[
    \Br(k_v) = \H^2(k_v, \Gbf_m) \longhookrightarrow \Qbf/\Zbf.
\]
For $\alpha \in \Br(k)$, we write $\inv_v(\alpha)$ for $\inv_v(\alpha|_{k_v})$.

\par\medskip
\noindent{\bf Varieties}\par

A variety over a field $K$ is a separated scheme of finite type. 
If $X/K$ is a smooth variety, a \emph{smooth compactification} of $X$ is a smooth proper variety containing $X$ as a dense open subset. By Nagata's compactification theorem and resolution of singularities, smooth compactifications exist; see \cite{Nagata63Resolution,Hironaka1964I,Hironaka1964II}. 

Denote by $\Z_0(X)$ the group of zero-cycles on $X$ ({\em i.e.} the free abelian group on the closed points of $X$). For $\delta \in \Zbf$, denote by $\Z_0^\delta(X) \subseteq \Z_0(X)$ the subset of zero-cycles of degree $\delta$. For a finite extension $L/K$, we write $\N_{L/K}$ for the pushforward of zero-cycles along $X_L\to X$. 

\par\medskip
\noindent{\bf Brauer groups}\par

By convention, the Brauer group of a scheme $X$ is the {\em cohomological Brauer group}
\[
    \Br(X) \coloneq \H_{\et}^2(X, \Gbf_m).
\]

For a smooth geometrically integral variety $X$ over a field $K$, denote by $\Br_{\nr}(X)$ the subgroup of $\Br(K(X))$ that consists of the {\em unramified} classes. This subgroup is also the Brauer group of any smooth compactification of $X$; {\em cf.} \cite[Proposition~6.2.7]{CTS2021Brauer}. Set
\[
    \Br_0(X) \coloneq \Img \bigl(\Br(K) \longrightarrow \Br(X)\bigr), \qquad \Br_1(X) \coloneq \Ker \bigl(\Br(X) \longrightarrow \Br(X_{\bar{K}}) \bigr),
\]
and 
\[
    \Br_{\nr,1}(X) \coloneq \Br_{\nr}(X) \cap \Br_1(X) = \Ker \bigl(\Br_{\nr}(X) \longrightarrow \Br_{\nr}(X_{\bar{K}})\bigr).
\]
We call $\Br_1(X)$ the algebraic Brauer group and
$\Br_{\nr,1}(X)$ the unramified algebraic Brauer group of $X$.

\par\medskip
\noindent{\bf Brauer--Manin pairing for zero-cycles}\par

Let $X$ be a smooth geometrically integral variety over a number field $k$. For each $v \in \Omega_k$, we have a local pairing
\[
    \langle -, -\rangle_v \colon \Z_0(X_{k_v}) \times \Br(X) \longrightarrow \Qbf/\Zbf, \qquad \left(\sum_i n_i P_i, \Acal \right) \longmapsto \sum_{i} n_i \inv_v \bigl(\Cor_{k_v(P_i)/k_v}(\Acal(P_i))\bigr),
\]
where $\Cor$ denotes corestriction. This pairing factors through rational equivalence \cite[Proposition~6.4.2]{CTS2021Brauer} if $X$ is proper. For a variety over a nonarchimedean local field $K$, we write $\langle-,-\rangle_K$ for the analogous pairing defined using $\inv_K$. For $X/k$ as above, restricting to unramified classes and summing over the places gives the Brauer--Manin pairing
\[
    \langle -,- \rangle \colon
    \prod_{v \in \Omega_k} \Z_{0}(X_{k_v}) \times \Br_{\nr}(X) \longrightarrow \Qbf/\Zbf, \qquad
    ((z_v)_v,\Acal) \longmapsto  \sum_{v \in \Omega_k} \langle z_v, \Acal \rangle_v,
\]
which vanishes on (the diagonal image of) $\Z_0(X)$; see \cite[\S{15.1}]{CTS2021Brauer}. For the ordinary Brauer--Manin pairing, every family of local points is orthogonal to $\Br_0(X)$, whereas for local zero-cycles this orthogonality is characterized by the following criterion.

\begin{lemm}\label{lemm:OrthogonalityToBr0}
   Let $z_\Omega = (z_v)_v$ be a family of local zero-cycles on $X$. Then $z_\Omega$ is orthogonal to $\Br_0(X)$ if and only if there exists $\delta \in \Zbf$ such that
    \[
        \deg(z_v) = \delta
    \]
    for every finite place $v$, and 
    \[
        \deg(z_v) \equiv \delta \pmod{2}
    \]
    for every real place $v$.
\end{lemm}
\begin{proof}
    Let $\pi \colon X \to \Spec(k)$ denote the structure morphism. Its pushforward is the degree map
\[
    \deg \colon \Z_0(X) \longrightarrow \Z_0(\Spec(k)) = \Zbf.
\]
Put $\delta_v\coloneq\deg(z_v)$. By \cite[(6.2)]{CTS2021Brauer}, for all $\alpha \in \Br(k)$, one has
\[
    \langle z_\Omega,\pi^\ast \alpha\rangle = \sum_{v \in \Omega_k} \langle z_v,\pi^\ast \alpha \rangle_v = \sum_{v \in \Omega_k} \langle 
    \deg(z_v), \alpha \rangle_v = \sum_{v \in \Omega_k} \delta_v \inv_v(\alpha) 
\]
in $\Qbf/\Zbf$. Assume the existence of  $\delta \in \Zbf$ as in the statement. Then $\delta_v \inv_v(\alpha) = \delta \inv_v(\alpha)$ for all $v \in \Omega_k$, because $\Br(\Rbf) = \Zbf/2\Zbf$ is $2$-torsion and $\Br(\Cbf) = 0$. By global reciprocity, one has
\[
    \langle z_\Omega, \pi^\ast \alpha\rangle = \delta \sum_{v \in \Omega_k} \inv_v(\alpha) = 0
\]
for all $\alpha \in \Br(k)$. Conversely, assume that $z_\Omega$ is orthogonal to $\Br_0(X)$. Fix distinct finite places $v_1,v_2 \in \Omega_k$. For any $n \ge 1$, the Albert--Brauer--Hasse--Noether exact sequence gives a class $\alpha \in \Br(k)$ with
\[
    \inv_v(\alpha) = \begin{cases}
         \frac{1}{n} & \text{if } v = v_1,\\
         -\frac{1}{n} & \text{if } v = v_2,\\
         0 & \text{otherwise}.
    \end{cases}
\]
Orthogonality to $\Br_0(X)$ then implies $\delta_{v_1} \equiv \delta_{v_2} \pmod{n}$ for all $n \ge 1$, and hence $\delta_{v_1} = \delta_{v_2}$. The integer $\delta = \delta_v$ is therefore independent of the finite place $v$. Similarly, taking $n=2$, $v_1$ any real place and $v_2$ any finite place, we have that $\delta_{v_1} \equiv \delta \pmod{2}$.
\end{proof}


\subsection{Unramified Brauer group} \label{subsec:BrauerGroup}

The following finiteness result for the unramified Brauer groups is standard. It is applicable to smooth projective compactifications of homogeneous spaces under connected linear algebraic groups. Indeed, connected linear algebraic groups are geometrically rational by Chevalley's theorem \cite{Chevalley1954Group}, so their homogeneous spaces are geometrically unirational, hence rationally connected, see \cite[\S~14.1, point~(3)]{CTS2021Brauer}. The following lemma follows from point~(6) of the same section. 

\begin{lemm}
\label{lemm:FinitenessBrNrRC}
The geometric Picard group $\Pic(Y_{\bar{k}})$ of any smooth projective rationally connected variety $Y$ over a field $k$ is finitely generated and torsion-free, and the quotient
\[
    \Br(Y)/\Br_0(Y)
\]
is finite.
\end{lemm}

\medskip
The following lemma shows how finite geometric stabilizers detect
geometric unramified Brauer classes on homogeneous spaces and will
be used in the proof of Proposition~\ref{prop:CorrectingPrimary}
to reduce to the case of a finite stabilizer.

\begin{lemm}\label{lemm:GeometricBrFiniteDectection}
    Let $\bar{K}$ be an algebraically closed field, $G$ a connected linear algebraic $\bar{K}$-group, and $H$ a closed subgroup of $G$. Then there exists a finite subgroup $F \subseteq H(\bar{K})$ such that the natural map
    \[
        \Br_{\nr}(G/H) \longrightarrow \Br_{\nr}(G/F)
    \]
    is injective. If $p$ is a prime number and $S \subseteq F$ is a $p$-Sylow subgroup, then the map
    \[
        \Br_{\nr}(G/H)\{p\} \longrightarrow \Br_{\nr}(G/S)\{p\}
    \]
    is injective, where $\{p\}$ denotes the $p$-primary subgroup.
\end{lemm}
\begin{proof}
    To prove the first assertion, it suffices to find a finite subgroup $F \subseteq H(\bar{K})$ such that $G/F \to G/H$ admits a rational section, or equivalently such that $(G/F)_\eta$ has an $L$-point, where $\eta $ is the generic point of $G/H$ and $L$ is the function field of $G/H$. In order to give an $L$-point of $(G/F)_\eta$, we consider the generic fibre of $G\to G/H$, namely
    \[
        X \coloneq G\times_{G/H}\Spec(L),
    \]
    which is an $L$-torsor under $H$ defining an element $[X]\in \H^1(L,H)$. It suffices to seek a finite subgroup $F\subseteq H(\bar{K})$ and an $L$-torsor $Y$ under $F$, together with an $F$-equivariant morphism $Y\to X$. Indeed, taking quotients by $F$ would then give the desired point
    \[
        \Spec(L) = Y/F \longrightarrow X/F \longrightarrow (G/F)_\eta.
    \]
    By a result of Chernousov--Gille--Reichstein \cite[Theorem~1.1]{CGR2006Torsor}, there exists a finite subgroup $F \subseteq H(\bar{K})$ such that the induced map
    \[
        \H^1(E,F) \longrightarrow \H^1(E,H)
    \]
    is surjective for every field extension $E/\bar{K}$. By our choice of $F$, there exists an $F_L$-torsor $Y$ together with an $F_L$-equivariant $L$-morphism $Y \to X$. As above, taking quotients yields an $L$-point of $(G/F)_\eta$.
    
    For the second assertion, if $S \subseteq F$ is a $p$-Sylow subgroup, then the finite \'etale morphism
    \[
        G/S \longrightarrow G/F
    \]
    has degree $[F:S]$, which is prime to $p$. Restriction--corestriction therefore gives injectivity on $p$-primary Brauer classes under pullback from $G/F$ to $G/S$; see \cite[Proposition~3.8.4]{CTS2021Brauer}. Composing this pullback with the injection in the first assertion shows that the map
    \[
        \Br_{\nr}(G/H)\{p\} \longrightarrow \Br_{\nr}(G/S)\{p\}
    \]
    is injective.
\end{proof}

\medskip
If $G = \GL_n$, Bogomolov's formula \cite{Bogomolov1988Brauer} identifies $\Br_{\nr}(G/F)$ with a subgroup of $\H^2(F,\Qbf/\Zbf)$. Let us show a similar result when $G$ has trivial Picard group.

\begin{prop}\label{prop:BrNrSubgroupOfH2}
    Let $\bar{K}$ be an algebraically closed field, $G$ a connected linear $\bar{K}$-algebraic group with $\Pic(G) = 0$, and let $F \subseteq G(\bar{K})$ be a finite subgroup. Then there is an injection
    \[
        j \colon \H^2(F, \bar{K}^\times) \longhookrightarrow \Br(G/F),
    \]
    where $\bar{K}^\times$ is endowed with the trivial $F$-action. Its image contains $\Br_{\nr}(G/F)$. Furthermore, if $F$ is a $p$-group for some prime $p$, then $\Br_{\nr}(G/F)$ is $p$-primary torsion.
\end{prop}
\begin{proof}
    Choose a closed embedding $G \hookrightarrow \GL(V)$ and take the closure of $G$ in $\Pbf(\End(V) \oplus \bar{K})$ via $g \mapsto [g:1]$. Right multiplication on the first summand extends the right translation action. A projective equivariant resolution that is an isomorphism over $G$ gives a smooth projective compactification $G \subseteq X$. Denote by $[X/F]$ the Deligne--Mumford quotient stack. By functoriality of the Hochschild--Serre spectral sequences associated with the $F$-torsors $X \to [X/F]$ and $G \to G/F$, the inclusion $G \hookrightarrow X$ induces a commutative diagram with exact rows
    \[
    \begin{tikzcd}[column sep=large, row sep=large]
        & \H^2(F,\bar{K}^\times) \arrow[r] \arrow[d] \arrow[dr, dashed, "j" description] & \Br([X/F]) \arrow[r] \arrow[d] & \H^1(F,\Pic(X)) \\
        0 \arrow[r] & \H^2(F,\bar{K}[G]^\times) \arrow[r] & \Br(G/F).
    \end{tikzcd}
    \]
    Here we use the assumption $\Pic(G) = 0$ and the fact that $\Br(X) = 0$ because $X$ is smooth, projective and rational by Chevalley's theorem \cite{Chevalley1954Group}.

    By birational invariance, $\Pic^0(X)=0$. Since $X$ is rational, Lemma~\ref{lemm:FinitenessBrNrRC} shows that $\Pic(X)$ is finitely generated and torsionfree. For any line bundle $L$ on $X$, the family $a^\ast L\otimes p_X^\ast L^{\vee}$ on $G\times X$, where $a$ denotes the action and $p_X$ the projection, induces a morphism
    \[
         \varphi_L \colon G \longrightarrow \Pic_{X /\bar{K}}, \qquad g \longmapsto [g^\ast L \otimes L^\vee];
    \]
    see \cite[\S{2.5}]{CTS2021Brauer}. Since $G$ is connected and $\varphi_L(e)=0$, this morphism factors through $\Pic^{0}_{X/\bar{K}}$. Thus
    \[
        [g^\ast L\otimes L^{\vee}] \in \Pic^{0}_{X/\bar{K}}(\bar{K}) = \Pic^{0}(X) = 0
    \]
    for every $g \in G(\bar{K})$. In other words, $F \subseteq G(\bar{K})$ acts trivially on $\Pic(X)$, so
    \[
        \H^1(F,\Pic(X)) = \Hom(F,\Pic(X)) = 0.
    \]
    The map $\H^2(F,\bar{K}^\times)\to\Br([X/F])$ is therefore surjective. By \cite[\S{4}]{KT2025Brauer}, we have inclusions
    \[
        \Br_{\nr}(G/F) \subseteq \Br([X/F]) \subseteq \Br(G/F).
    \]
    Consequently, the induced map $\H^2(F, \bar{K}^\times) \to \Br(G/F)$ has image containing $\Br_{\nr}(G/F)$. By the diagram, its injectivity follows from that of
    \[
        \H^2(F,\bar{K}^\times) \longrightarrow \H^2(F,\bar{K}[G]^\times),
    \]
    which we now prove.
    
    By Rosenlicht's lemma \cite[Proposition~3]{Rosenlicht1957Rationality}, we have a short exact sequence
    \[
        1 \longrightarrow \bar{K}^\times \longrightarrow \bar{K}[G]^\times \longrightarrow \Xbb^\ast(G) \longrightarrow 0,
    \]
    where $\Xbb^\ast(G)$ denotes the group of characters $G \longrightarrow \Gbf_m$.
    Note that the translation action of $F$ induces the trivial action on $\Xbb^\ast(G)$, because
    \[
        (f^\ast\chi)(g) = \chi(gf) = \chi(g)\chi(f)
    \]
    for all $\chi \in \Xbb^\ast(G)$, $f \in F$, and $g \in G(\bar{K})$, and hence $f^\ast \chi$ and $\chi$ are equal modulo a constant. Since $G$ is connected, $\Xbb^\ast(G)$ is a finitely generated torsion-free abelian group. Hence\[
        \H^1(F,\Xbb^\ast(G)) = \Hom(F,\Xbb^\ast(G)) = 0,
    \]
    which implies the injectivity of
    \[
        \H^2(F,\bar{K}^\times) \longrightarrow \H^2(F,\bar{K}[G]^\times).
    \]
    This establishes the first assertion of the proposition. 
    
    Finally, the group $\H^2(F,\bar{K}^\times)$ is $|F|$-torsion, hence $p$-primary if $F$ is a $p$-group.
\end{proof}

\medskip
Let $K$ be a field of characteristic zero, $\bar{K}$ an algebraic closure of $K$, and $G$ a connected
linear algebraic $K$-group. Let $F\subseteq G(K)$ be a finite constant $K$-subgroup, and put $Y \coloneq G/F$. Recall the canonical isomorphism
\[
    \H^2(F,K^\times) \simeq \Ext^{c}_{K\gp}(F,\Gbf_m),
\]
where $F$ acts trivially on $K^\times$ and where
$\Ext^{c}_{K\gp}(F,\Gbf_m)$ denotes
the group of equivalence classes of {\em central} extensions of $F$ by $\Gbf_m$ as algebraic $K$-groups; see \cite[proof of Proposition~1.1]{Dem17}.
The construction of Colliot-Th\'el\`ene--Xu
\cite[\S~2, pp.~313--314]{CTX09} associates with the
torsor $G\to Y$ a homomorphism
\[
    \rho_{\tors(G)} \colon \Ext^{c}_{K\gp}(F, \Gbf_m)
    \longrightarrow \Br(Y),
\]
which takes the class of a central extension
\[
    e \colon \qquad 1 \longrightarrow \Gbf_m \longrightarrow E \longrightarrow F \longrightarrow 1
\]
to the boundary $\partial_e([G \to Y])$, where $\partial_e \colon \H^1_{\et}(Y, F) \to \H^2_{\et}(Y, \Gbf_m)$ is the associated connecting map.
This construction is compatible with base change.
Lemma~\ref{lemm:BrauerCentralExtensionCompatibility} identifies it, over $\bar{K}$, with the homomorphism of Proposition~\ref{prop:BrNrSubgroupOfH2}.

\medskip
We work in the setting of Proposition~\ref{prop:BrNrSubgroupOfH2},
with $G$ and $F$ defined over $\bar K$. Because $\bar{K}$ is algebraically closed, one has a noncanonical isomorphism
\[
     \Qbf/\Zbf \cong \Qbf/\Zbf(1) \coloneq \varinjlim_n \mu_n \subseteq \bar{K}^\times.
\]
Choose such an isomorphism; the resulting identifications of cohomology groups over this algebraic closure use this choice.
Furthermore, the quotient $\frac{\bar{K}^\times}{\Qbf/\Zbf(1)}$ is uniquely divisible, hence {\em cohomologically trivial}. Thus
\[
   \H^2(F,\Qbf/\Zbf) \cong \H^2(F,\Qbf/\Zbf(1)) \cong \H^2(F,\bar{K}^\times).
\]
Proposition~\ref{prop:BrNrSubgroupOfH2} then says that, provided $\Pic(G) = 0$, every element $\Acal \in \Br_{\nr}(G/F)$ can be represented by a central extension
\[
    0 \longrightarrow \Qbf/\Zbf \longrightarrow F_\Acal \longrightarrow F \longrightarrow 1.    
\]

In finite group cohomology, the coefficients $\Zbf$ and $\Qbf/\Zbf$ carry the trivial action. For further computational purposes, it is useful to exploit the following description of $\H^2(F,\Qbf/\Zbf)$ using the {\em nonabelian exterior square} of $F$. Recall that $F \wedge F$ is the group generated by symbols $x\wedge y$, for $x,y\in F$, subject to the following ``crossed bilinearity'' and alternating relations
\begin{align*}
    (xx')\wedge y & =({}^{x}x'\wedge{}^{x}y)(x\wedge y),\\
    x\wedge(yy') &=(x\wedge y)({}^{y}x\wedge{}^{y}y'),\\
    x\wedge x &=1,
\end{align*}
for all $x,x',y,y'\in F$, where ${}^{x}y \coloneq xyx^{-1}$. The formula
\[
    \kappa\colon F\wedge F\longrightarrow [F,F], \qquad x\wedge y\longmapsto [x,y] \coloneq xyx^{-1}y^{-1}
\]
yields a surjective homomorphism. Its kernel
\[
    \M(F)\coloneq\Ker(\kappa)
\]
is canonically isomorphic to the {\em Schur multiplier} $\H_2(F,\Zbf)$; see \cite[\S~2]{Moravec2012Brauer}. On the other hand, we have
\[
    \H^2(F,\Qbf/\Zbf) \simeq \Hom(\M(F),\Qbf/\Zbf)
\]
by the universal coefficient theorem. The following lemma gives an explicit isomorphism between these groups.

\begin{lemm}\label{lemm:MoravecFormula}
   Let $F$ be a finite group and let $\Acal \in \H^2(F,\Qbf/\Zbf)$ be represented by a central extension $0 \to \Qbf/\Zbf \to F_\Acal \to F \to 1$. Then there is an isomorphism 
   \[
    \H^2(F,\Qbf/\Zbf) \overset{\cong}{\longrightarrow} \Hom(\M(F),\Qbf/\Zbf)
    \]
    given by 
    \[
        \Acal \longmapsto \Bigl(
        m = (a_1\wedge b_1) \cdots (a_s\wedge b_s)
        \longmapsto [\widetilde a_1,\widetilde b_1] \cdots [\widetilde a_s,\widetilde b_s] \Bigr),
    \]
    where $m\in\M(F)$, $a_1,\ldots,a_s,b_1,\ldots,b_s \in F$ and $\widetilde{a}_1,\ldots,\widetilde{a}_s, \widetilde{b}_1,\ldots,\widetilde{b}_s \in F_{\Acal}$ are their respective liftings. 
\end{lemm}
\begin{proof}
    This isomorphism is constructed in \cite[proof of Theorem~3.2]{Moravec2012Brauer}. 
\end{proof}

\section{Existence of zero-cycles on homogeneous spaces} \label{sec:ExistenceOfZeroCycles}

This section is devoted to the proof of Theorem~\ref{thm:Main}.


\subsection{A local flexibility lemma} \label{subsec:LocalFlexibility}

The aim of this subsection is to prove
Proposition~\ref{prop:LocalFlexibility}.
We first give an explicit formula for the local invariant.

Fix a prime number $p$. Throughout this subsection, let $K$ be a $p$-adic local field
such that
\begin{equation}\label{eq:LocalRootsOfUnity}
    \begin{cases}
        \mu_p\subseteq K & \text{if } p\ne2,\\
        \mu_4\subseteq K & \text{if } p=2.
    \end{cases}
\end{equation}
Let $K(p)$ be its maximal $p$-extension and put
$\Gamma\coloneq\Gal(K(p)/K)$. Let $q$ be the largest power of $p$ such that $\mu_q\subseteq K$. Fix isomorphisms $\mu_q \cong \Zbf/q\Zbf \cong \frac{1}{q}\Zbf/\Zbf$ of $\Gamma$-modules. Throughout this subsection, we use these isomorphisms to identify both the coefficient groups and the corresponding cohomology groups. Since $\mu_q$ is $p$-primary, the inflation $\H^2(\Gamma,\mu_q) \to \H^2(K,\mu_q)$ is an isomorphism by \cite[Exercise~8.2(d)]{Harari2017}. Via the Kummer identification, we have an isomorphism
\[
	\H^2(\Gamma, \mu_{q}) = \H^2(K,\mu_{q}) \xrightarrow{\inv_K} \Zbf/q\Zbf \subseteq \Qbf/\Zbf
\] 
by local class field theory. The next lemma unwinds the description of $\inv_K$ in terms of central extensions of $\Gamma$. Under our assumptions on $K$, Labute's classification of Demu\v{s}kin groups \cite[Theorem~7]{Labute67Demushkin} gives the pro-$p$ presentation

\begin{equation} \label{eq:GaloisGroupOfMaximalPExtension}
	\Gamma = \langle x_1,\ldots,x_{d+2} \mid x_1^q[x_1^{-1},x_2^{-1}]\cdots[x_{d+1}^{-1},x_{d+2}^{-1}] = 1 \rangle,
\end{equation}
where $d \coloneq [K:\Qbf_p]$ is even.

\begin{lemm}\label{lemm:DescriptionOfLocalInvariant}
	There exists an integer $h$ prime to $p$ such that the above local invariant map $\inv_K\colon \H^2(\Gamma, \mu_{q}) \longrightarrow \Zbf/q\Zbf$ is described as follows. Let $\alpha \in \H^2(\Gamma,\mu_q)$ be represented by a central extension
	\[
		1 \longrightarrow \mu_q \longrightarrow \Gamma_\alpha \longrightarrow \Gamma \longrightarrow 1
	\]
	of pro-$p$ groups, and choose arbitrary lifts $\widetilde{x}_i \in \Gamma_\alpha$ of $x_i$ (for $1 \le i \le d+2$). Then 
    \[
        \inv_K(\alpha) = h \cdot \widetilde{x}_1^{\,q} [\widetilde{x}_1^{-1},\widetilde{x}_2^{-1}] \cdots [\widetilde{x}_{d+1}^{-1},\widetilde{x}_{d+2}^{-1}] \quad \text{in } \mu_q \simeq \Zbf/q\Zbf.
    \]
\end{lemm}
\begin{proof}
	In view of \eqref{eq:GaloisGroupOfMaximalPExtension}, the pro-$p$ group $\Gamma$ fits into an exact sequence
	\begin{equation} \label{eq:DescriptionOfLocalInvariant:1}
		1 \longrightarrow \Rcal \longrightarrow \Fcal \longrightarrow \Gamma \longrightarrow 1,
	\end{equation}
	where $\Fcal$ is the free pro-$p$ group on $d+2$ generators $X_1,\ldots,X_{d+2}$, and $\Rcal \subseteq \Fcal$ is the closed normal subgroup topologically generated by the relator $R \coloneq X_1^q[X_1^{-1},X_2^{-1}]\cdots[X_{d+1}^{-1},X_{d+2}^{-1}]$. The Hochschild--Serre spectral sequence associated with \eqref{eq:DescriptionOfLocalInvariant:1} yields a bijective transgression map 
	\[
		\tg\colon \H^1(\Rcal,\Zbf/q\Zbf)^{\Fcal} \overset{\cong}{\longrightarrow}\H^2(\Gamma,\Zbf/q\Zbf);
	\]
	see \cite[\S{1.4}]{Labute67Demushkin}. Here the group $\H^1(\Rcal,\Zbf/q\Zbf)^{\Fcal}$ consists of the continuous $\Fcal$-equivariant homomorphisms $\Rcal \to \Zbf/q\Zbf$. The composition of the inverse transgression with evaluation at $R^{-1}$
	\[
		\bar{R} \colon \H^2(\Gamma,\Zbf/q\Zbf) \longrightarrow \Zbf/q\Zbf, \qquad \alpha \longmapsto \tg^{-1}(\alpha)(R^{-1})
	\]
	is a group isomorphism; {\em cf.} \cite[p.~115]{Labute67Demushkin}. Since $\inv_K$ and $\bar R$ are both group isomorphisms onto $\Zbf/q\Zbf$, there exists an integer $h$ prime to $p$ such that 
    \[
        \inv_K(\alpha) = h \cdot \bar{R}(\alpha)
    \]
    for all $\alpha \in \H^2(\Gamma,\Zbf/q\Zbf) \cong \H^2(\Gamma,\mu_q)$. 
    \par
    
    Since $\Fcal$ is free, the lifts $\widetilde{x}_i$ define a continuous homomorphism
    \[
        \widetilde{\pi}\colon \Fcal \longrightarrow \Gamma_\alpha
    \]
    lifting the projection $\pi \colon \Fcal \to \Gamma$. By the group cohomology result in Lemma~\ref{lemm:TransgressionLiftedRelator} below, we have $\bar{R}(\alpha) = \widetilde{\pi}(R)$. It follows that
    \[
	    \inv_K(\alpha)
        = h \cdot \bar{R}(\alpha) = h \cdot \widetilde{\pi}(R) = h \cdot \widetilde{x}_1^q [\widetilde{x}_1^{-1},\widetilde{x}_2^{-1}]
        \cdots [\widetilde{x}_{d+1}^{-1},\widetilde{x}_{d+2}^{-1}]
    \]
	in $\mu_q \cong \Zbf/q\Zbf$, which proves the lemma.
\end{proof}

\begin{lemm}\label{lemm:TransgressionLiftedRelator}
    Let $1 \to \Rcal \to \Fcal \overset{\pi}{\rightarrow} \Gamma \to 1$ be an exact sequence of pro-$p$ groups, let $A$ be a finite abelian group on which $\Gamma$ acts trivially, and assume that the transgression map
    \[
        \tg\colon \H^1(\Rcal,A)^{\Fcal} \overset{\cong}{\longrightarrow}\H^2(\Gamma,A)
    \]
    is bijective. Let $\alpha\in\H^2(\Gamma,A)$ be represented by a central extension $1 \to A \to \Gamma_\alpha \to \Gamma \to 1$, and let $\widetilde{\pi}\colon\Fcal\to\Gamma_\alpha$ be a continuous homomorphism lifting $\pi$. For every $R\in\Rcal$, one has
    \[
        \tg^{-1}(\alpha)(R^{-1})=\widetilde{\pi}(R)\in A.
    \]
\end{lemm}
\begin{proof}
    By \cite[Proposition~1.6.6]{NSW13}, it suffices to construct a continuous $1$-cochain $u' \colon \Fcal \to A$ whose coboundary is the inflation of a cocycle representing $\alpha$ and whose restriction to $\Rcal$ is given by $R \mapsto \widetilde\pi(R)^{-1}$.

    Choose a continuous section $\sigma \colon \Gamma \to \Gamma_\alpha$ with $\sigma_1=1$, and set
    \[
        u'_W \coloneq \sigma_{\pi(W)} \widetilde{\pi}(W)^{-1} \in A.
    \]
    Since $\widetilde\pi$ is a homomorphism and $A$ is central, $du'$ is the inflation of the cocycle $\sigma_s\sigma_t\sigma_{st}^{-1}$ representing $\alpha$. For every $R\in\Rcal$, we have $\pi(R)=1$, so $u'_R = \widetilde\pi(R)^{-1}$ as required.
\end{proof}

\begin{prop}[Local flexibility] \label{prop:LocalFlexibility}
    Let $p$ be a prime and let $K$ be a $p$-adic local field satisfying \eqref{eq:LocalRootsOfUnity}. Let $q$ be the largest power of $p$ such that $\mu_q \subseteq K$. We use the identification $\mu_q \simeq \frac{1}{q}\Zbf/\Zbf \subseteq \Qbf/\Zbf$ fixed above. Let $S$ be a finite $p$-group, $B \subseteq \H^2(S,\Qbf/\Zbf)$ a subgroup, and
    \[
        \theta \colon B \longrightarrow \Qbf/\Zbf
    \]
    a functional. Assume that $K$ is sufficiently large in the following sense:
    \begin{enumerate}
        \item \label{prop:LocalFlexibility:1} The exponents of both $S$ and $B$ divide $q$.
            
        \item \label{prop:LocalFlexibility:2} $d = [K:\Qbf_p]$ is sufficiently large (larger than a constant depending only on $S$).
    \end{enumerate}
    Then there exists a continuous homomorphism 
    \[
        \rho \colon \Gal(\bar{K}/K) \longrightarrow S
    \] 
    such that for every $\alpha \in B$ that comes from an element $\widetilde{\alpha} \in \H^2(S,\mu_q) \cong \H^2(S,\Zbf/q\Zbf)$, one has
    \[
        \inv_K(\rho^\ast \widetilde{\alpha}) = \theta(\alpha).
    \]
\end{prop}

\noindent Note that every element of $B$ is $q$-torsion and hence comes from $\H^2(S,\Zbf/q\Zbf)$.

\begin{proof}
    Let $h$ be the integer from Lemma~\ref{lemm:DescriptionOfLocalInvariant}. Since $\Qbf/\Zbf$ is a divisible abelian group, the homomorphism $h^{-1}\cdot\theta$ extends to some homomorphism
    \[
        \theta' \colon \H^2(S,\Qbf/\Zbf) \longrightarrow \Qbf/\Zbf,
    \]
    where $h^{-1}$ denotes the inverse of $h$ modulo $q$. Since $\M(S)=\H_2(S,\Zbf)$ is finite, there exists an integer $s_0$, depending only on $S$, such that every element of $\M(S)$ can be written as
    \[
        (a_1 \wedge b_1) \cdots(a_s \wedge b_s), \qquad s \le s_0, \qquad a_i,b_i \in S.
    \]
    By Pontryagin duality, $\theta'$ is evaluation at some
    \[
        m = (a_1 \wedge b_1)\cdots(a_s \wedge b_s)\in \M(S)
    \]
    with $s\le s_0$ and $[a_1,b_1] \cdots [a_s,b_s]=1$. Assume $d\ge 2s_0$. We may define a homomorphism 
    \begin{align*}
        \rho\colon \Gamma \longrightarrow S, \qquad & \rho(x_3) = a_1^{-1}, \ \rho(x_4) = b_1^{-1},\ \ldots,\ \rho(x_{2s+1}) = a_s^{-1},\ \rho(x_{2s+2}) = b_s^{-1}, \\
        & \rho(x_1) = \rho(x_2) = 1, \qquad \rho(x_i) = 1 \quad (2s+3\le i\le d+2),
    \end{align*}
    in view of \eqref{eq:GaloisGroupOfMaximalPExtension}.
    Let $\widetilde{\alpha} \in \H^2(S,\mu_q) \cong \H^2(S,\Zbf/q\Zbf)$, and let $\alpha$ be its image in $\H^2(S,\Qbf/\Zbf)$. Represent $\widetilde{\alpha}$ by a central extension
    \[
        1 \longrightarrow \mu_{q} \longrightarrow S_{\widetilde{\alpha}} \longrightarrow S \longrightarrow 1.
    \]
    Pulling this extension back along $\rho$ yields a central extension
    \[
    	1 \longrightarrow \mu_{q} \longrightarrow \Gamma_{\rho^\ast \widetilde{\alpha}} \longrightarrow \Gamma \longrightarrow 1
    \]
    representing $\rho^\ast \widetilde{\alpha} \in \H^2(\Gamma,\mu_q)$. For $1 \le i \le d+2$, lift the generator $x_i \in \Gamma$ to $\widetilde{x}_i \in \Gamma_{\rho^\ast \widetilde{\alpha}}$. Then
    \[
        \inv_K(\rho^\ast \widetilde{\alpha}) = h \cdot \widetilde{x}_1^q [\widetilde{x}_1^{-1},\widetilde{x}_2^{-1}] \cdots [\widetilde{x}_{d+1}^{-1},\widetilde{x}_{d+2}^{-1}] \in \mu_q \cong \Zbf/q\Zbf \subseteq \Qbf/\Zbf
    \]
    by Lemma~\ref{lemm:DescriptionOfLocalInvariant}. The pullback extension also comes with a homomorphism
    \[
        \Gamma_{\rho^\ast\widetilde{\alpha}} \longrightarrow S_{\widetilde{\alpha}}.
    \]
    For $i=1,2$ or $2s+3\le i\le d+2$, the image of $\widetilde{x}_i$ in $S_{\widetilde{\alpha}}$ lies in the central subgroup $\mu_q$ because $\rho(x_i)=1$. Let $\widetilde{a}_1, \widetilde{b}_1,\ldots,\widetilde{a}_s, \widetilde{b}_s \in S_{\widetilde{\alpha}}$ be the respective images of $\widetilde{x}_3^{-1},\widetilde{x}_4^{-1},\ldots,\widetilde{x}_{2s+1}^{-1},\widetilde{x}_{2s+2}^{-1}$. These are respective lifts of $a_1,b_1,\ldots,a_s,b_s$. After extending coefficients along $\mu_q \hookrightarrow \Qbf/\Zbf$, Lemma~\ref{lemm:MoravecFormula} gives
    \[
        h \cdot \theta'(\alpha) = h \cdot  \alpha(m) =  h \cdot  [\widetilde{a}_1,\widetilde{b}_1] \cdots [\widetilde{a}_s,\widetilde{b}_s] = h \cdot \widetilde{x}_1^q [\widetilde{x}_1^{-1},\widetilde{x}_2^{-1}] \cdots  [\widetilde{x}_{d+1}^{-1},\widetilde{x}_{d+2}^{-1}] =\inv_K(\rho^\ast \widetilde{\alpha})
    \]
    in $\mu_q \cong \Zbf/q\Zbf \subseteq \Qbf/\Zbf$. In particular, if $\alpha \in B$ then
    \[
        \theta(\alpha) = \inv_K(\rho^\ast \widetilde{\alpha}).
    \]
    Composing $\rho$ with the quotient map $\Gal(\bar{K}/K)\twoheadrightarrow\Gamma$ completes the proof.
\end{proof}


\subsection{Correcting primary Brauer classes by local closed points} \label{subsec:CorrectingPrimary}

The main result of this subsection is
Proposition~\ref{prop:CorrectingPrimary}.
We first treat the case of a finite stabilizer over a $p$-adic field (Proposition~\ref{prop:CorrectingPrimaryFinite}), using the local flexibility Proposition~\ref{prop:LocalFlexibility}.
We then reduce the general case to this setting, using the central cover provided by Lemma~\ref{lemm:SpecialCover} and the finite subgroup provided by Lemma~\ref{lemm:GeometricBrFiniteDectection}.

Let $K$ be a $p$-adic field with algebraic closure $\bar{K}$, and let $\widetilde{G}$ be a connected linear algebraic $K$-group. Let $S \subseteq \widetilde{G}(K)$ be a finite subgroup, viewed as a constant $K$-group. Put $n \coloneq |S|$, $Y \coloneq \widetilde{G}/S$. Let $B \subseteq \Br_{\nr}(Y)$ be a finite subset. All extensions lie inside $\bar{K}$. For finite $M/K$ containing $\mu_n$, write
$\Gamma_M \coloneq \Gal(\bar{K}/M)$. For a continuous homomorphism $\rho \colon \Gamma_M \to S$ and $a \in \H^2(S,\mu_n(\bar{K}))$, write $\rho^\ast a$ for its image under
\[
    \H^2(S,\mu_n(\bar{K})) \overset{\rho^\ast}{\longrightarrow} \H^2(M,\mu_n) \longhookrightarrow \Br(M),
\]
where $S$ acts trivially on $\mu_n(\bar{K})$. 

\begin{lemm}
\label{lemm:FlexibilityZeroCycles}
    With the above notation, assume that $\Pic(\widetilde{G}_{\bar{K}}) = 0$. Let $\bar{B}$ be the geometric image of $B$. For each $\bar{\Acal} \in \bar{B}$, choose $a_{\bar{\Acal}} \in \H^2(S,\mu_n(\bar{K}))$ mapping to $\bar{\Acal}$ under
    \[
        \H^2(S,\mu_n(\bar{K})) \longrightarrow \H^2(S,\bar{K}^\times) \stackrel{j}{\longrightarrow} \Br(Y_{\bar{K}}),
    \]
    where $j$ is the map of Proposition~\ref{prop:BrNrSubgroupOfH2} and the $S$-actions on the coefficients are trivial. Then there exists a finite extension $L/K$ containing $\mu_n$ such that the following assertions hold for every finite extension $M/L$.
    \begin{enumerate}
        \item \label{lemm:FlexibilityZeroCycles:1} For every $Q \in Y(M)$, let $\rho_Q \colon \Gamma_M \to S$ be a continuous homomorphism representing the class of the fibre $\widetilde{G}_M\times_{Y_M}Q$ in $\H^1(M,S)$. Then
        \[
            \Acal_M(Q) = \rho_Q^\ast a_{\bar{\Acal}} \quad \text{in } \Br(M) \qquad (\Acal \in B).
        \]
        
        \item \label{lemm:FlexibilityZeroCycles:2} If moreover $\H^1(M,\widetilde{G})=1$, then for every continuous homomorphism $\rho \colon \Gamma_M \to S$, there exists a zero-cycle $c \in \Z_0(Y)$ of degree zero such that
        \[
            \langle c,\Acal \rangle_K = \inv_M(\rho^\ast a_{\bar{\Acal}}) \qquad (\Acal\in B).
         \]
    \end{enumerate}
\end{lemm}
\begin{proof}
    Put $K_0 \coloneq K(\mu_n)$. We construct an auxiliary class $\beta_{\bar{\Acal}}$ to compute the evaluation of $\Acal$.

    \medskip
    \noindent
    \emph{Step 1. Constructing Brauer classes from the chosen lifts.}

    The lifts $a_{\bar\Acal}$ in the statement exist by Proposition~\ref{prop:BrNrSubgroupOfH2} and the fact that $n \H^2(S,\bar{K}^\times) = 0$. We use the notation and construction recalled in Subsection~\ref{subsec:BrauerGroup}. For each $\bar{\Acal}\in\bar{B}$, let
    \[
        e_{\bar {\Acal}} \colon \qquad 1 \longrightarrow \Gbf_m \longrightarrow E_{\bar {\Acal}} \longrightarrow S_{K_0} \longrightarrow 1
    \]
    be a central extension of $K_0$-group schemes corresponding to the image of $a_{\bar{\Acal}}$ under
    \[
        \H^2(S,\mu_n(\bar{K})) \longrightarrow \H^2(S,K_0^\times) \simeq \Ext^{c}_{K_0\gp}(S_{K_0},\Gbf_m).
    \]
    Set
    \[
        \beta_{\bar {\Acal}} \coloneq \rho_{\tors(\widetilde{G}_{K_0})}([e_{\bar {\Acal}}]) = \partial_{e_{\bar {\Acal}}}([\widetilde{G}_{K_0} \longrightarrow Y_{K_0}])\in\Br(Y_{K_0}).
    \]

    \medskip 
    \noindent 
    \emph{Step 2. Identifying these classes with the original classes after a finite extension.}

    By compatibility with base change and Lemma~\ref{lemm:BrauerCentralExtensionCompatibility}, we have
    \[
        \beta_{\bar {\Acal}}|_{\bar{K}} = \rho_{\tors(\widetilde{G}_{\bar{K}})}
    ([(e_{\bar {\Acal}})_{\bar{K}}]) = \bar{\Acal}.
    \]
    Consequently,
    \[
        \Acal_{K_0} - \beta_{\bar {\Acal}} \in \Ker\!\left(\Br(Y_{K_0}) \longrightarrow \Br(Y_{\bar{K}}) \right).
    \]
    By \cite[Proposition~6.6.10]{Poonen2017Rational} and the finiteness of $B$, there exists a common finite extension $L/K_0$ inside $\bar{K}$ such that $\Acal_L = \beta_{\bar {\Acal}}|_L \in\Br(Y_L)$ for every $\Acal\in B$. In particular, for every finite extension $M/L$ and every $Q \in Y(M)$,
    \[
        \Acal_M(Q) = (\beta_{\bar{\Acal}})_M(Q)
    \qquad (\Acal\in B).
    \] 
    We now aim to identify the right-hand side with the boundary $\partial_{(e_{\bar{\Acal}})_M}([\rho_Q])$ and show that this boundary equals $\rho_Q^\ast a_{\bar{\Acal}}$.

    \medskip
    \noindent
    \emph{Step 3. Computing evaluation at a point.}

    Fix a finite extension $M/L$. The torsor $\widetilde{G}_M \to Y_M$ defines the map
    \[
        \ev_{\widetilde{G}_M} \colon Y(M) \longrightarrow \H^1(M,S), \qquad Q \longmapsto [\widetilde{G}_M \times_{Y_M} Q].
    \]
    The first two rows of the diagram in \cite[Proposition~2.7]{CTX09} give the commutative diagram
    \[
    \begin{array}{ccccc}
        Y(M) & \times & \Br(Y_M) & \xrightarrow{\mathrm{evaluation}} & \Br(M) \\[3mm]
        {\scriptstyle \ev_{\widetilde{G}_M}} \big \downarrow && \big\uparrow{\scriptstyle \rho_{\tors(\widetilde{G}_M)}} && \Vert \\[3mm] \H^1(M,S) & \times & \Ext^{c}_{M\text{-gp}}(S_M,\Gbf_m) & \xrightarrow{\mathrm{boundary}} & \Br(M).
    \end{array}
    \]
    The lower horizontal pairing takes $(\xi,[e])$ to $\partial_e(\xi)$. Let $Q\in Y(M)$, and choose a continuous homomorphism $\rho_Q \colon \Gamma_M \to S$ representing the class of its fibre. Since the upward arrow takes $[(e_{\bar {\Acal}})_M]$ to $(\beta_{\bar {\Acal}})_M$, the diagram yields
    \[
        (\beta_{\bar {\Acal}})_M(Q) = \partial_{(e_{\bar {\Acal}})_M}([\rho_Q])
        \quad \text{in } \Br(M).
    \]
    The extension $(e_{\bar{\Acal}})_M$ corresponds to the image of $a_{\bar{\Acal}}$ under
    \[
        \H^2(S, \mu_n(\bar{K})) \longrightarrow \H^2(S, M^\times) \overset{\cong}{\longrightarrow} \Ext^{c}_{M\gp}(S_M, \Gbf_m).
    \]
    For any central extension $e$ of $S_M$ by $\Gbf_m$ with middle term $E$, Hilbert's Theorem~90 gives $\H^1(M, \Gbf_m)=0$, so $E(M) \to S$ is surjective. Taking $M$-points of $(e_{\bar{\Acal}})_M$ gives the central extension of abstract groups $1 \rightarrow M^\times \rightarrow E_{\bar{\Acal}}(M) \rightarrow S \rightarrow 1$. Then we have
    \[
        \H^2(S,M^\times) \overset{\cong}{\longrightarrow} \Ext^{c}_{M\gp}(S_M,\Gbf_m) \xrightarrow{[e]\mapsto[e(M)]} \Ext^{c}_{\mathrm{grp}}(S,M^\times).
    \]
    Choose a normalized cocycle $c \colon S \times S \to \mu_n(\bar{K})$ representing $a_{\bar{\Acal}}$. Both the $M$-group extension $(e_{\bar{\Acal}})_M$ and the extension of abstract groups obtained by taking its $M$-points are represented by this cocycle via the following multiplication on $\Gbf_m \times S_M$
    \[
        (u,g)(v,h)=(uv\,c(g,h),gh).
    \]
    With this multiplication, $s(g) = (1,g)$ defines a section $s \colon S \to E_{\bar{\Acal}}(M)$. Since each $s(g)$ is an $M$-point of $E_{\bar{\Acal}}$, it is fixed by $\Gamma_M$. Hence the boundary $\partial_{(e_{\bar{\Acal}})_M}([\rho_Q])$ is represented by
    \[
        (\sigma,\tau)\longmapsto s(\rho_Q(\sigma))\,\sigma\bigl(s(\rho_Q(\tau))\bigr)\,s(\rho_Q(\sigma\tau))^{-1}=c\bigl(\rho_Q(\sigma),\rho_Q(\tau)\bigr).
    \]
    Hence $\partial_{(e_{\bar{\Acal}})_M}([\rho_Q]) = \rho_Q^\ast a_{\bar{\Acal}}$ in $\Br(M)$. Combining the equalities above gives $\Acal_M(Q) = \rho_Q^\ast a_{\bar{\Acal}}$ for every $\Acal\in B$, which proves ~\ref{lemm:FlexibilityZeroCycles:1}.

    \medskip
    \noindent
    Now let $\rho \colon \Gamma_M \to S$ be the given homomorphism. It defines a class $[\rho] \in \H^1(M,S)$. The exact sequence of pointed sets
    \[
        \widetilde{G}(M) \longrightarrow Y(M)
        \overset{\ev_{\widetilde{G}_M}}{\longrightarrow}
        \H^1(M, S) \longrightarrow \H^1(M, \widetilde{G})
    \]
    and the assumption $\H^1(M, \widetilde{G}) = 1$ provide a point $Q \in Y(M)$ whose fibre represents $[\rho]$. Thus $\rho_Q$ is conjugate to $\rho$. Conjugation does not change the pullback of group cohomology classes with trivial coefficients. Let $Q_0 \in Y(K)$ be the image of the identity. Its fibre over $M$ is the trivial $S_M$-torsor, so~\ref{lemm:FlexibilityZeroCycles:1} gives $\Acal_M(Q_0) = 0$ for every $\Acal \in B$. Hence
    \[
        \Acal_M(Q) - \Acal_M(Q_0) = \rho^\ast a_{\bar{\Acal}} \quad \text{in } \Br(M) \qquad(\Acal\in B).
    \]
    For the projection $f \colon Y_ M\to Y$, set
    \[
        c \coloneq f_\ast\bigl([Q] - [Q_{0,M}]\bigr) \in \Z_0(Y),
    \]
    which has degree zero. The composite
    \[
        \Br(M) \overset{\Cor_{M/K}}{\longrightarrow} \Br(K) \overset{\inv_K}{\longrightarrow} \Qbf/\Zbf
    \]
    equals $\inv_M$. For every $\Acal \in B$, we therefore have
    \begin{align*}
     \langle c, \Acal \rangle_K & = \inv_K\! \left(\Cor_{M/K} \bigl(\Acal_M(Q) - \Acal_M(Q_0)\bigr)\right) \\
        & = \inv_M\bigl(\Acal_M(Q) - \Acal_M(Q_0)\bigr) \\
        & = \inv_M(\rho^\ast a_{\bar{\Acal}}).
    \end{align*}
    This proves~\ref{lemm:FlexibilityZeroCycles:2}.
\end{proof}

\begin{prop}[Correction for finite stabilizers]
\label{prop:CorrectingPrimaryFinite}
    Let $K$ be a $p$-adic local field and let $\widetilde{G}$ be a connected linear algebraic group over $K$ such that
    \[
        \Pic(\widetilde{G}_L) = 0 \qquad \text{and} \qquad \H^1(L,\widetilde{G})=1
    \]
    for every finite extension $L/K$. Let $F \subseteq \widetilde{G}$ be a finite $K$-subgroup and put $X \coloneq \widetilde{G}/F$. Let $B \subseteq \Br_{\nr}(X)$ be a finite $p$-primary subgroup and let
    \[
        \theta \colon \bar{B} \longrightarrow \Qbf/\Zbf
    \]
    be a homomorphism, where $\bar{B}$ denotes the geometric image of $B$.

    Then there exists a zero-cycle $c \in \Z_0(X)$ of degree zero such that 
    \[
        \langle c,\Acal\rangle_K = \theta(\bar \Acal)
        \qquad \text{for every } \Acal \in B,
    \]
    where $\bar \Acal$ denotes the geometric image of $\Acal$.
\end{prop}

\begin{proof}
    In order to apply Proposition~\ref{prop:LocalFlexibility}, let $S$ be a Sylow $p$-subgroup of $F(\bar{K})$. Take a finite extension $K'/K$ inside $\bar{K}$ splitting $F$, view $S$ as a constant subgroup of $F_{K'}$, and put $Y \coloneq \widetilde{G}_{K'}/S$. The inclusion $S \subseteq F_{K'}$ induces a natural morphism $\pi \colon Y \to X_{K'}$.

    Since $\pi_{\bar{K}}$ has degree $[F_{\bar{K}}:S]$ prime to $p$, \cite[Proposition~3.8.4]{CTS2021Brauer} implies that $\Ker(\pi_{\bar{K}}^\ast)$ has trivial $p$-primary part. As $\bar{B}$ is $p$-primary, the restriction of $\pi_{\bar{K}}^\ast$ to $\bar{B}$ is injective. Together with Proposition~\ref{prop:BrNrSubgroupOfH2}, this gives the following injections:
    \[
    \begin{tikzcd}[column sep=large]
        \bar{B}
        \arrow[r, hook, "\pi_{\bar{K}}^\ast"']
        \arrow[rr, bend left=25, "\iota"]
        & \Br_{\nr}(Y_{\bar{K}})
        \arrow[r, hook] &
    \H^2(S,\Qbf/\Zbf)
    \end{tikzcd}
    \]
    Here $\iota$ denotes the composite injection, via which we identify $\bar{B}$ with a subgroup of $\H^2(S,\Qbf/\Zbf)$. Set $B_Y \coloneq \pi^\ast(B|_{K'}) \subset\Br_{\nr}(Y)$ and $n \coloneq |S|$. For each geometric class $\bar {D} \in\overline{B_Y}$, choose a lift
    \[
        a_{\bar{D}}\in \H^2(S,\mu_n(\bar{K}))
    \]
    of $\bar{D}$, as in Lemma~\ref{lemm:FlexibilityZeroCycles}. Applying that lemma over $K'$ to $Y$, $B_Y$, and these lifts gives a finite extension $L/K'$, containing $K'(\mu_n)$, over which its conclusions hold.

    Choose a finite extension $M/L$ satisfying \eqref{eq:LocalRootsOfUnity} with $K = M$ and conditions~\ref{prop:LocalFlexibility:1} and~\ref{prop:LocalFlexibility:2} of Proposition~\ref{prop:LocalFlexibility}. Let $q$ be the largest power of $p$ such that $\mu_q\subset M$. Since $\mu_n\subset L$, we have $n\mid q$. We also write $a_{\bar D}$ for its image in $\H^2(S,\mu_q(\bar{K}))$. Using for $\mu_q$ the restriction of the coefficient identification used to define $\iota$, we apply Proposition~\ref{prop:LocalFlexibility} over $M$ to $S$, $\iota(\bar{B})$, and $\theta \circ \iota^{-1} \colon \iota(\bar{B}) \to \Qbf/\Zbf$. This gives a continuous homomorphism
    \[
        \rho \colon \Gal(\bar{K}/M) \longrightarrow S
    \]
    such that
    \[
        \inv_M \left(\rho^\ast a_{\pi_{\bar{K}}^\ast(\bar{\Acal})}\right) = (\theta \circ \iota^{-1}) \bigl(\iota(\bar{\Acal})\bigr) = \theta(\bar{\Acal}) \qquad\text{for every }\Acal\in B.
    \]

    Since $\H^1(M, \widetilde{G}) = 1$, Lemma~\ref{lemm:FlexibilityZeroCycles}~\ref{lemm:FlexibilityZeroCycles:2}, applied over $K'$ with the extension $M/L$ to $B_Y$, gives a zero-cycle $c_Y \in \Z_0(Y)$ of degree zero such that
    \[
        \left\langle c_Y, \pi^\ast(\Acal|_{K'})\right\rangle_{K'} = \inv_M\left(\rho^\ast a_{\pi_{\bar{K}}^\ast(\bar\Acal)}\right) = \theta(\bar\Acal)
    \]
    for every $\Acal\in B$.

    Set
    \[
        c \coloneq \N_{K'/K}(\pi_\ast c_Y) \in \Z_0(X).
    \]
    Then $\deg(c)=0$, and by functoriality of the Brauer pairing,
    \[
        \langle c, \Acal\rangle_K = \langle c_Y, \pi^\ast(\Acal|_{K'}) \rangle_{K'} = \theta(\bar\Acal).
    \]
    This proves the proposition.
\end{proof}

\begin{lemm} \label{lemm:SpecialCover}
    Let $K$ be a $p$-adic local field and let $G$ be a {\em split} connected linear algebraic group over $K$. Then there exists a
    central isogeny
    \[
        \widetilde{G} \longrightarrow G,
    \]
    where $\widetilde{G}$ is a connected linear algebraic group over $K$ such that
    \[
        \Pic\bigl(\widetilde{G}_{K'}\bigr) = 0 \qquad \text{and} \qquad \H^1\bigl(K',\widetilde{G}_{K'}\bigr) = 1
    \]
    for any finite extension $K'/K$.
\end{lemm}
\begin{proof}
    Let $G = U \rtimes L$ be a Levi decomposition, where $U$ is unipotent and $L$ is split reductive; such a decomposition exists by Mostow's theorem \cite{Mostow1956Group}. The radical $T\coloneq\Rad(L)$ is a central torus in $L$, hence is contained in any maximal torus of $L$. Since $L$ contains a split maximal torus by assumption, $T$ is a split torus as well. On the other hand, the derived subgroup $\Der(L)$ is semisimple and $L$ is an almost direct product of $\Rad(L)$ and $\Der(L)$, {\em i.e.}, the multiplication
    \[
        \Der(L) \times T \longrightarrow L, \qquad (s,t) \longmapsto st
    \]
    is a central isogeny. Let $\widetilde{L}$ be the simply connected cover of $\Der(L)$ and put 
    \[
        \widetilde{G}\coloneq U \rtimes \bigl(\widetilde{L} \times T\bigr).
    \]
    We obtain a central isogeny $\widetilde{G} \to G$. 
    
    Now, we have $\Pic(U_{K'}) = 0$ because $U$ is an affine space. Second, $\Pic\bigl(\widetilde{L}_{K'}\bigr) = 0$ by \cite[Lemme~6.9~(iv)]{Sansuc1981Arithmetique}. Finally, $\Pic(T_{K'}) = 0$ because $T$ is split and $\Pic(\Gbf_m) = 0$. It follows that $\Pic\bigl(\widetilde{G}_{K'}\bigr) = 0$, by Sansuc's units--Picard--Brauer exact sequence \cite[(6.10.1)]{Sansuc1981Arithmetique}. 

    On the other hand, $\H^1(K',U) = 1$ by \cite[Lemme~1.13]{Sansuc1981Arithmetique}, $\H^1\bigl(K', \widetilde{L}\bigr) = 1$ by a classical result of Kneser--Bruhat--Tits \cite{Kneser1965I,Kneser1965II,BT1987Group}, and $\H^1(K',T) = 0$ by Hilbert's Theorem 90. Thus $\H^1\bigl(K', \widetilde{G}\bigr) = 1$ by the long exact sequence for nonabelian Galois cohomology.
\end{proof}

\begin{prop}[Correction in the general case]
\label{prop:CorrectingPrimary}
    Let $K$ be a $p$-adic local field and let $X$ be a homogeneous space of a connected linear algebraic group $G$ over $K$. Let $B \subseteq \Br_{\nr}(X)$ be a finite $p$-primary subgroup, and denote by $\bar{B}$ its geometric image. Let
    \[
        \theta \colon \bar{B} \longrightarrow \Qbf/\Zbf
    \]
    be a homomorphism. Then there exists a zero-cycle $c \in \Z_0(X)$ of degree zero such that
    \[
        \langle c, \Acal \rangle_K = \theta(\bar{\Acal}) \qquad \text{for every } \Acal \in B,
    \]
    where $\bar \Acal$ denotes the geometric image of $\Acal$.
\end{prop}

\begin{proof}
    Take a finite extension $K'/K$ inside $\bar{K}$ such that $X(K') \neq \varnothing$ and $G_{K'}$ is split. Let $P^0\in X(K')$. By Lemma~\ref{lemm:SpecialCover}, after replacing $G_{K'}$ by a central cover, we obtain a connected linear algebraic group $\widetilde{G}$ such that $\Pic(\widetilde{G}_L)=0$ and $\H^1(L,\widetilde{G})=1$ for every finite extension $L/K'$.

    Let $\widetilde{H} \subseteq \widetilde{G}$ be the stabilizer of $P^0$. By Lemma~\ref{lemm:GeometricBrFiniteDectection}, there exists a finite subgroup $\bar{F} \subseteq \widetilde{H}(\bar{K})$ such that
    \[
        \Br_{\nr}(X_{\bar{K}})\{p\} \longrightarrow \Br_{\nr}(\widetilde{G}_{\bar{K}}/\bar F)\{p\}
    \]
    is injective. After enlarging $K'$ if necessary, we may assume that $\bar F$ comes from a finite $K'$-subgroup $F \subseteq \widetilde{H}$.

    Put $Y \coloneq \widetilde{G}/F$, and let $\pi \colon Y \longrightarrow X_{K'}$ be the natural morphism. Set
    \[
        B_Y \coloneq \pi^\ast(B|_{K'})\subseteq\Br_{\nr}(Y).
    \]
    Since $B$ is $p$-primary and $\pi_{\bar{K}}^\ast$ is injective on the $p$-primary part, it induces an isomorphism
    \[
        \pi_{\bar{K}}^\ast \colon \bar{B} \xrightarrow{\cong} \overline{B_Y}.
    \]
    Let $\theta_Y\colon \overline{B_Y}\to \mathbf Q/\mathbf Z$ be the functional corresponding to $\theta$ under this isomorphism.

    Let $c_Y\in \Z_0(Y)$ be a zero-cycle of degree zero provided by Proposition~\ref{prop:CorrectingPrimaryFinite} applied to $Y$, $B_Y$, and $\theta_Y$. Set $c \coloneq \N_{K'/K}(\pi_\ast c_Y)\in \Z_0(X)$. Then $c$ has degree zero. For every $\Acal \in B$, the projection formula and the identity $\inv_K \circ \Cor_{K'/K} = \inv_{K'}$ give
    \begin{align*}
        \langle c,\Acal\rangle_K & =\langle \pi_\ast c_Y,\Acal|_{K'}\rangle_{K'}\\
        & = \langle c_Y,\pi^\ast(\Acal|_{K'})\rangle_{K'}\\
        & = \theta_Y\bigl(\overline{\pi^\ast(\Acal|_{K'})}\bigr) \\
        & = \theta_Y\bigl(\pi_{\bar{K}}^\ast(\bar \Acal)\bigr) \\
        & = \theta(\bar \Acal).
    \end{align*}
    This proves the proposition.
\end{proof}


\subsection{Proof of the main theorem} \label{subsec:GlobalCorrection}

We are now ready to prove our main result. 

\begin{proof}[Proof of Theorem~\ref{thm:Main}]
    We are given a homogeneous space $X$ of a connected linear algebraic group over a number field $k$, an integer $\delta \in \Zbf$, and a family
    \[
        z_\Omega \in \prod_{v \in \Omega_k} \Z_0^\delta(X_{k_v})
    \]
    of local zero-cycles that is orthogonal to $\Br_{\nr,1}(X)$. By Lemma~\ref{lemm:FinitenessBrNrRC}, there exists a finite subgroup $B \subseteq \Br_{\nr}(X)$ that generates $\Br_{\nr}(X)$ modulo $\Br_0(X)$. Let 
    \[
        \bar{B}\coloneq \frac{B}{B \cap \Br_{\nr,1}(X)} \subseteq \Br_{\nr}(X_{\bar{k}})
    \]
    denote its geometric image. To apply the local correction result of Proposition~\ref{prop:CorrectingPrimary}, we treat the primary components of $B$ separately. For each prime $p \mid |B|$, write $|B| = p^a m$ with $\gcd(p,m)=1$, and use the Chinese remainder theorem to choose an integer $e_p$ such that
    \[
        e_p \equiv 1 \pmod{p^a} \qquad \text{and} \qquad e_p \equiv 0 \pmod{m}.
    \]
    Then $e_p \Acal = \Acal_p$ is the $p$-primary component of $\Acal$ for every $\Acal\in B$. For each such $p$, orthogonality of $z_\Omega$ to $\Br_{\nr,1}(X)$ implies that the homomorphism
    \[
        B\{p\} \longrightarrow \Qbf/\Zbf, \qquad \Acal_p \longmapsto \langle z_\Omega,\Acal_p \rangle
    \]
    factors through a homomorphism $\theta_p \colon \bar{B}\{p\} \to \Qbf/\Zbf$. Choose a place $v_p \in \Omega_k$ lying above $p$ and apply Proposition~\ref{prop:CorrectingPrimary} to $X_{k_{v_p}}$, the image of $B\{p\}$ in $\Br_{\nr}(X_{k_{v_p}})$, and $\theta_p$. We obtain a zero-cycle $c_p \in \Z_0(X_{k_{v_p}})$ of degree zero such that
    \[
        \langle c_p,\Acal_p \rangle_{v_p} = \theta_p(\bar \Acal_p) \qquad (\Acal_p \in B\{p\}).
    \]
    Set $c'_p \coloneq e_p c_p$. Then for every $\Acal \in B$,
    \[
        \langle c'_p,\Acal\rangle_{v_p} = \langle c_p,\Acal_p\rangle_{v_p} = \theta_p(\bar \Acal_p) = \langle z_\Omega,\Acal_p\rangle.
    \]
    For $v \in \Omega_k$, set
    \[
        z'_v \coloneq \begin{cases}
            z_v - c'_p & \text{if } v=v_p \text{ for some prime } p \mid |B|,\\
            z_v & \text{otherwise}.
        \end{cases}
    \]
Since $\deg(c'_p) = 0$ for every $p$, the family
    \[
        z'_\Omega \coloneq (z'_v)_v \in \prod_{v \in \Omega_k} \Z_0^\delta(X_{k_v})
    \]
    still consists of zero-cycles of degree $\delta$. For any $\Acal \in B$, we have
    \[
        \langle z'_\Omega,\Acal \rangle = \langle z_\Omega,\Acal \rangle - \sum_{p \mid |B|}\langle c'_p,\Acal \rangle_{v_p} = 0.
    \]
    Thus the family $z'_\Omega$ is orthogonal to $B$. Since this family is also orthogonal to $\Br_0(X)$ by Lemma~\ref{lemm:OrthogonalityToBr0}, it is orthogonal to $\Br_{\nr}(X)$. It follows from \cite[Th\'eor\`eme~A]{HW2020Galois} that $\Z_0^\delta(X^c) \neq \varnothing$ for any smooth compactification $X^c$ of $X$. By Lemma~\ref{lemm:ChowMovingLemma} below, one has $\Z_0^\delta(X) \neq \varnothing$.
\end{proof}

\begin{lemm}[Chow's moving lemma for zero-cycles]
\label{lemm:ChowMovingLemma}
    Let $Y$ be a smooth variety over a field $K$, and let $U\subseteq Y$ be a dense open subset. For every zero-cycle $z \in \Z_0(Y)$, there exists a zero-cycle $z'$ supported in $U$ that is rationally equivalent to $z$. In particular, if $Y$ is proper, then
    \[
        \deg(z')=\deg(z).
    \]
\end{lemm}
\begin{proof}
    This follows from \cite[Complement to \S{3}, p.~599]{CT2005Finitude}.
\end{proof}

\section{Transcendental Brauer--Manin obstruction revisited} \label{sec:Example}


\subsection{Optimality of the algebraic Brauer group} \label{subsec:Optimality}

Even when the unramified algebraic Brauer group of a torsor under a torus has order two modulo constants, the assumption about orthogonality to $\Br_{\nr,1}$ in Theorem~\ref{thm:Main} cannot be relaxed further.

\begin{exam}[Colliot-Thélène] \label{exam:Optimality}
    We recall Colliot-Thélène's example in \cite{CT2014Brauer}. Consider the principal homogeneous space $Y$ under a five-dimensional torus over $\Qbf$. The space $Y$ is given in \cite[Remarque~5.2]{CT2014Brauer} by the equation
    \[
        (x^2-13y^2)(z^2-17t^2)(u^2-221w^2)=5.
    \]
    By \cite[Theorem~4.1 and Proposition~5.1~(a)]{CT2014Brauer}, we have $\prod_{v \in \Omega_\Qbf} Y(\Qbf_v) \neq \varnothing$ and
    \[
        \Br_{\nr,1}(Y) / \Br_0(Y) \cong \Zbf/2\Zbf.
    \]
    This quotient is generated by the class of the Azumaya algebra $\Acal = (x^2-13y^2,17)$.
    Since $17$ is a square in $\bar{\Qbf}$, the class of $\Acal$ vanishes in $\Br(Y_{\bar{\Qbf}})$ and is therefore algebraic.
    \medskip

    We now show that $\Acal$ gives a Brauer--Manin obstruction to the existence of zero-cycles of odd degree. At every place $v \in \Omega_\Qbf$, either $13$, $17$, or $221$ is a square in $\Qbf_v$. Therefore \cite[Proof of Proposition~5.1~(d)]{CT2014Brauer} shows that $\Acal_v$ is constant: more precisely, it is the pullback of
    \[
        \beta_v \coloneq \begin{cases}
            (5,17)_v \in \Br(\Qbf_v) & \text{if } 13 \in \Qbf_v^{\times2},\\
            0 & \text{otherwise}.
        \end{cases}
    \]
    For every local zero-cycle $z_v \in \Z_0(Y_{\Qbf_v})$, restriction--corestriction therefore gives
    \[
        \langle z_v, \Acal \rangle_v = \deg(z_v) \inv_v(\beta_v).
    \]
    Since $(5,17)$ has zero invariant at every place other than $5$ and $17$, it suffices to compute $\inv_v(\beta_v)$ at these two places. Since $13$ is a square in $\Qbf_{17}$ but not in $\Qbf_5$ and $\left(\frac{5}{17}\right) = -1$, we obtain
    \[
        \inv_v(\beta_v) = \begin{cases}
            \frac{1}{2} & \text{if } v=17,\\
            0 & \text{otherwise}.
        \end{cases}
    \]
    Consequently, for every $\delta\in\Zbf$, every family $z_\Omega \in \prod_{v \in \Omega_{\Qbf}} \Z_0^{\delta}(Y_{\Qbf_v})$ satisfies
    \[
        \langle z_\Omega, \Acal \rangle = \frac{\delta}{2} \in \Qbf/\Zbf.
    \]
    Thus no such family of odd degree is orthogonal to $\Br_{\nr,1}(Y)$, so the global reciprocity law excludes global zero-cycles of odd degree.
\end{exam}

Therefore, in Example~\ref{exam:Optimality}, the group $\Br_{\nr,1}(Y)$ is optimal for Theorem~\ref{thm:Main} in the sense that it cannot be replaced by any proper subgroup containing $\Br_0(Y)$.


\subsection{Failure of the Hasse principle with perfect geometric stabilizers} \label{subsec:ExamplePerfectStabilizers}

A cheap way to construct homogeneous spaces with an algebraic Brauer group (not only the unramified algebraic part) reduced to constants is to take finite {\em perfect} stabilizers.

\begin{coro}
\label{coro:HassePerfectStabilizer}
    Let $k$ be a number field, let $\delta \in \Zbf$, and let $X$ be a homogeneous space of a simply connected semisimple linear algebraic $k$-group whose geometric stabilizers $\bar{H}$ are finite and perfect. Then
    \[
        \prod_{v\in\Omega_k} \Z_0^\delta(X_{k_v}) \neq \varnothing \Longrightarrow \Z_0^\delta(X) \neq \varnothing.
    \]  
\end{coro}
\begin{proof}
    By assumption we have $\bar{H}^{\ab} = 1$. Therefore $\Br_{1}(X) = \Br_0(X)$ by \cite[Théorème~1]{Dem10}, and {\em a fortiori} $\Br_{\nr,1}(X) = \Br_0(X)$. By Lemma~\ref{lemm:OrthogonalityToBr0}, every equidegree family of local zero-cycles is orthogonal to $\Br_0(X)$. The conclusion follows from Theorem~\ref{thm:Main}.
\end{proof}




The cup-product construction from \cite{Linh2026Insufficiency} can be generalized to obtain homogeneous spaces with perfect geometric stabilizers that are counterexamples to the Hasse principle. By Corollary~\ref{coro:HassePerfectStabilizer}, these spaces possess zero-cycles of degree one. The following lemma transforms this construction into a study of group cohomology.

\begin{lemm}\label{lemm:HomoSpacesCoho}
    Let $k$ be a field of characteristic zero and $Q$ a finite $k$-group scheme with center $Z$, and let $\mathbf n_Q \in \H^2(k,Q)$ be its distinguished neutral class. For $\beta\in\H^2(k,Z)$, let $X_\beta$ be a homogeneous space of $\SL_N$ with Springer class $\beta \cdot \mathbf n_Q$, obtained from the Demarche--Lucchini Arteche construction. For every extension $K/k$ we have
    \[
        X_\beta(K)\neq\varnothing \quad \Longleftrightarrow \quad \beta|_K \in \Img\bigl(\Delta_K \colon \H^1(K,Q/Z) \longrightarrow \H^2(K,Z)\bigr).
    \]
    Here $\Delta_K$ is the connecting map associated with the central extension $1 \to Z \to Q \to Q/Z \to 1$.
\end{lemm}
\begin{proof}
    Since $\H^1(K,\SL_N)=1$, the space $X_\beta$ has a $K$-point if and only if its Springer class $(\beta\cdot\mathbf n_Q)|_K$ is neutral. By \cite[Proposition~2.3]{Borovoi1993Abelianization}, this is equivalent to $\beta|_K\in\Img(\Delta_K)$.
\end{proof}

Let $G$ be a finite group, $q$ a power of a prime $p$, and $V$ a representation of $G$ of dimension $n > 0$ over $\Fbf_q$. Assume that either $Z(G) = 1$ or the action of $G$ on $V$ is faithful. This assumption will be used in Lemma~\ref{lemm:GroupQ}. Fix an integer $d\geq2$, to be chosen sufficiently large; see
Proposition~\ref{prop:NoRationalPoint} and
Subsection~\ref{subsec:ExampleNilpotentStabilizers} for the respective bounds. Set
\[
    M \coloneq V^d \qquad \text{and} \qquad M^\vee \coloneq (V^\vee)^d,
\]
with the diagonal $G$-actions, where $V^\vee \coloneq \Hom_{\Fbf_q}(V,\Fbf_q)$ and where $G$ acts on $V^\vee$ by 
\[
    (g\cdot\lambda)(v) \coloneq \lambda(g^{-1}v).
\]
All tensor products of $\Fbf_q$-modules below are taken over $\Fbf_q$. Let 
\[
    \pi \colon \Mat_d(\Fbf_q) \longrightarrow Z \coloneq \Mat_d(\Fbf_q)/\Fbf_q \cdot I_d
\]
be the quotient map and set $\phi \coloneq \pi \circ \psi$, where $\psi$ is defined by
\[
    \psi \colon M \otimes M^\vee \longrightarrow \Mat_d(\Fbf_q), \qquad x \otimes y \longmapsto (y_i(x_j))_{1\leq i,j\leq d},
\]
where $x = (x_1,\ldots,x_d)$ and $y =(y_1,\ldots,y_d)$, and where $Z$ carries the trivial $G$-action. Let $F \coloneq Z \times (M\oplus M^\vee)$ with the group law
\[
    (z, x, y)\cdot(z', x', y') = \bigl( z+z'+\phi(x\otimes y'), x+x', y+y'\bigr).
\]
Then $G$ acts on $F$ by group automorphisms via $g \cdot (z,x,y) \coloneq (z,gx,gy)$. 

Let $k$ be a number field with $\mu_p \subseteq k$ and fix an isomorphism $\Zbf/p\Zbf \simeq \mu_p$ of Galois modules. We regard $Q \coloneq F \rtimes G$ as a finite constant $k$-group. A class in $\H^1(k,(M \oplus M^\vee)\rtimes G)$ is represented by a cocycle $(x,y,g)$, where $g \colon \Gamma_k \to G$ is a homomorphism and
\[
    x \in \Z^1(k,{}_gM), \qquad y \in \Z^1(k,{}_gM^\vee).
\]
Here ${}_gM$ and ${}_gM^\vee$ denote the modules endowed with the {\em twisted} Galois action $\sigma \cdot m \coloneq g_\sigma m$. The cup product $[x] \cup [y]$ takes values in $\H^2(k,{}_gM\otimes{}_gM^\vee)$, and $\phi_\ast$ denotes the map on cohomology induced by $\phi$. 

We now construct $\beta\in \H^2(k, Z)$. We choose $d$ distinct finite places $v_1,\ldots,v_d$, none of which lies over $p$. For $1\leq i\leq d$, let $E_{i,i}\in\Mat_d(\Fbf_q)$ denote the matrix whose $(i,i)$-entry is $1$ and whose other entries are all $0$. For a matrix $E$, write $[E]\coloneq\pi(E)\in Z$. Let
\[
    ([E_{v}])_v\in \bigoplus_{\substack{v\text{ finite} \\ \text{or real}}} Z
\]
be the family with $[E_v] = [E_{i,i}]$ when $v=v_i$ and $0$ at all other places. Since $\sum_v[E_{v}]=[I_d]=0$ in $Z$, Lemma~\ref{lemm:GlobalReciprocity} gives a class $\beta\in\H^2(k,Z)$ such that $\operatorname{Inv}_v^Z(\beta|_{k_v})=[E_v]$ for every finite or real place $v$. Let $X=X_\beta$ be the corresponding homogeneous space from Lemma~\ref{lemm:HomoSpacesCoho}.

To calculate the center and the derived subgroup of $Q$, we use the following elementary facts about this construction.
\begin{enumerate}
\renewcommand{\labelenumi}{\textup{(\alph{enumi})}}
\renewcommand{\theenumi}{\alph{enumi}}
    \item \label{fact:GroupF} $Z(F)=Z$ and $[F,F]=Z$. Indeed, the map $\phi$ is {\em nondegenerate} in the sense of \cite[\S{3B}]{Linh2024BK} and surjective. Its nondegeneracy gives $Z(F)=Z$. Since $F/Z$ is abelian and
    \[
        [(0,x,0),(0,0,y)]=(\phi(x\otimes y),0,0),
    \]
    its surjectivity gives $[F,F]=Z$.
    
    \item \label{fact:NoInnerAut} Let $\rho \colon G \to \Aut(F)$ denote the action. Then $\rho(G) \cap \Inn(F) = \{\id_F\}$ because inner automorphisms act trivially on $F/Z = M \oplus M^\vee$, and any element of $G$ acting trivially on this quotient acts trivially on $F$.
\end{enumerate}

\begin{lemm}\label{lemm:GroupQ}
    With the above notation, the following are true.
    \begin{enumerate}
        \item \label{lemm:GroupQ:1} $Z(Q) = Z$.
        \item \label{lemm:GroupQ:2} If $G$ is perfect and $V$ is irreducible of dimension at least two (in particular, $G \neq 1$) then $Q$ is perfect.
    \end{enumerate}
    \end{lemm}
\begin{proof}
    We prove~\ref{lemm:GroupQ:1}. By fact~\eqref{fact:NoInnerAut} we have
    \[
        Z(Q) = Z(F)^G \times \bigl(Z(G) \cap \ker(\rho)\bigr).
    \]
    By assumption, either $Z(G)=1$ or the action of $G$ on $V$ is faithful, so the second factor is trivial. Hence $Z(Q) = Z$ by fact~\eqref{fact:GroupF}.

    We prove~\ref{lemm:GroupQ:2}. Assume $G$ perfect and $V$ irreducible and nontrivial. Since $Q = F \rtimes G$, we have $Q^{\ab} \simeq (F^{\ab})_G\oplus G^{\ab}$. Here the subscript $G$ means the coinvariant quotient. By fact~\eqref{fact:GroupF} we have $F^{\ab} = M \oplus M^\vee$, whose coinvariants vanish because $V$ and $V^\vee$ are nontrivial irreducible representations. Also $G^{\ab} = 0$ by assumption, so $Q^{\ab} = 0$, i.e. $Q$ is perfect.
\end{proof}

\begin{lemm}\label{lemm:ConnectingMap}
    With the above notation, we can therefore consider the connecting map associated with the central extension $1 \to Z \to Q \to (M \oplus M^\vee) \rtimes G \to 1$, which is given by 
    \begin{align*}
        \Delta \colon \H^1(k, (M\oplus M^\vee) \rtimes G) & \longrightarrow \H^2(k, Z)\\
        [(x,y,g)] & \longmapsto \phi_\ast([x] \cup [y]).
    \end{align*}
\end{lemm}
\begin{proof}
    Let $a = (x,y,g)$ be a cocycle representing a class in $\H^1(k, (M\oplus M^\vee) \rtimes G)$, and choose the lifts $\widetilde{a}_\sigma \coloneq (0,x_\sigma,y_\sigma,g_\sigma)\in Q$. The cocycle conditions are
    \[
        x_{\sigma\tau} = x_\sigma+g_\sigma x_\tau, \qquad y_{\sigma\tau} = y_\sigma+g_\sigma y_\tau, \qquad g_{\sigma\tau}=g_\sigma g_\tau.
    \]
    Hence the multiplication law in $Q$ yields
    \[
        \widetilde{a}_\sigma \widetilde{a}_\tau = \bigl(\phi(x_\sigma\otimes g_\sigma y_\tau), x_{\sigma\tau}, y_{\sigma\tau}, g_{\sigma\tau}\bigr).
    \]
    Since $Q$ is constant, the class $\Delta([a])$ is represented by the $Z$-valued cocycle
    \[
        (\sigma,\tau)\longmapsto \widetilde{a}_\sigma \widetilde{a}_\tau \widetilde{a}_{\sigma\tau}^{-1} = \phi(x_\sigma\otimes g_\sigma y_\tau),
    \]
    where we identify $Z$ with its image in $Q$. This cocycle is precisely the image under $\phi$ of the cocycle representing $[x] \cup [y]$, since the Galois action on ${}_gM^\vee$ is given by $\sigma\cdot y=g_\sigma y$. Thus $\Delta([a]) = \phi_\ast([x]\cup[y])$. The same formula holds over every extension of $k$.
\end{proof}

For each class in $\H^1(k,G)$, choose a representing cocycle $g$ and define
\[
    \Delta_g \colon \H^1(k,{}_gM\oplus{}_gM^\vee) \longrightarrow \H^2(k,Z), \qquad [(x,y)]\longmapsto \phi_\ast([x]\cup [y]).
\]
With this notation, the image of $\Delta$ is given by
\[
    \Img(\Delta) = \bigcup_{[g] \in \H^1(k, G)} \Img(\Delta_g).
\]
In view of Lemma~\ref{lemm:HomoSpacesCoho} and Lemma~\ref{lemm:ConnectingMap}, we now prove that the homogeneous space $X$ constructed above does not satisfy the Hasse principle.
\begin{lemm}
    For every place $v$ of $k$, there exist $x_v\in\H^1(k_v,M)$ and $y_v\in\H^1(k_v,M^\vee)$ such that $\beta|_{k_v}=\phi_\ast(x_v\cup y_v)$. Consequently, $X(k_v)\neq\varnothing$ for every $v$.
\end{lemm}
\begin{proof}
We adapt the proof of \cite[Lemma~5]{Linh2026Insufficiency}. If $v \notin \{v_1,\ldots,v_d\}$, take $x_v = y_v = 0$. For $v = v_i$, choose $a, b \in \H^1(k_v,\Fbf_p)$ with
$\Inv_v^{\Fbf_p}(a\cup b) = 1$ by Tate's local duality \cite[Théorème~10.9]{Harari2017}, and $w \in V$, $\ell \in V^\vee$ with $\ell(w) = 1$. Let $x_v$ and $y_v$ have respective $i$-th components
\[
    (x_v)_i=[\,\sigma\mapsto a_\sigma w\,], \qquad (y_v)_i=[\,\sigma\mapsto b_\sigma\ell\,],
\]
all other components being zero.
Let $j \colon \Fbf_p \hookrightarrow \Fbf_q$ be the natural inclusion. Then $\psi_\ast(x_v\cup y_v)=j_\ast (a\cup b) E_{i,i}$, so by functoriality of the local isomorphism in Lemma~\ref{lemm:GlobalReciprocity} we have
\[
    \Inv_v^{\Mat_d(\Fbf_q)} \bigl(\psi_\ast(x_v \cup y_v)\bigr) = \Inv_v^{\Fbf_q}\bigl(j_\ast(a \cup b)\bigr)E_{i,i} = E_{i,i}.
\]
Hence $\phi_\ast (x_v\cup y_v)=\beta|_{k_v}$.
Taking the trivial $G$-cocycle and applying lemmas~\ref{lemm:ConnectingMap}
and~\ref{lemm:HomoSpacesCoho}, we obtain $X(k_v)\neq\varnothing$.
\end{proof}

\begin{prop}\label{prop:NoRationalPoint}
    Assume $d > \max\limits_{v \mid p} n([k_v:\Qbf_p] + 2)$. Then for all $g \in \H^1(k, G)$, $x \in \H^1(k, {}_gM)$, and $y \in \H^1(k,{}_gM^\vee)$, we have $\beta \neq \phi_\ast(x \cup y)$. Consequently, $X(k)=\varnothing$.
\end{prop}
\begin{proof}
    We adapt the proof of \cite[Proposition~6]{Linh2026Insufficiency}. Suppose $\beta = \phi_\ast(x\cup y)$, where $x \in \H^1(k,{}_gM)$ and $y \in \H^1(k,{}_gM^\vee)$. For each place $v$, let $x_v \coloneq x|_{k_v}$ and $y_v \coloneq y|_{k_v}$. Then $\beta|_{k_v}=\pi_\ast(\psi_\ast(x_v\cup y_v))$.  By functoriality of the local isomorphism in Lemma~\ref{lemm:GlobalReciprocity}, we have the following commutative diagram for every finite or real place, where $\Inv_v^{\Mat}$ abbreviates $\Inv_v^{\Mat_d(\Fbf_q)}$. 
    \[
    \begin{tikzcd}
        \H^2(k_v,{}_gM\otimes {}_gM^\vee) \arrow[r, "\psi_\ast"] & \H^2(k_v,\Mat_d(\Fbf_q)) \arrow[r, "\pi_\ast"] \arrow[d, "\operatorname{Inv}^{\Mat}_v", "\simeq"'] & \H^2(k_v,Z) \arrow[d, "\operatorname{Inv}^{Z}_v", "\simeq"'] \\
        & \Mat_d(\Fbf_q) \arrow[r, "\pi"] & Z
    \end{tikzcd}
    \]
    For each finite or real place $v$, set $A_v \coloneq \Inv^{\Mat}_v(\psi_\ast(x_v \cup y_v))$. If $v = v_i$ for some $i$, then $\Inv^{Z}_v(\beta|_{k_v}) = [E_{i,i}]$ by construction, so the definition of $\pi$ gives 
    \begin{equation}\label{eq:A_v}
        A_v = \lambda_v I_d + E_{i,i}, \qquad \lambda_v \in \Fbf_q.
    \end{equation}
    If $v\notin\{v_1,\ldots,v_d\}$, we have
    \begin{equation} \label{eq:A_vOtherPlaces}
        A_v=\lambda_v I_d,\qquad \lambda_v\in\Fbf_q 
    \end{equation}
    instead. Now we study the rank of $A_v$. The class $x_v \in \H^1(k_v,{}_gM)$ determines a linear map $x_v \colon \Fbf_q^d \longrightarrow \H^1(k_v,{}_gV)$ via the identifications
    \[
        \H^1(k_v,{}_gM) \simeq (\Fbf_q^d)^\vee\otimes_{\Fbf_q}\H^1(k_v,{}_gV) \simeq\Hom_{\Fbf_q} \bigl(\Fbf_q^d,\H^1(k_v,{}_gV)\bigr).
    \]
    Here the first identification uses the standard basis for $\Fbf_q^d$. Similarly, Tate's local duality gives the isomorphism $\H^1(k_v,{}_gV^\vee) \simeq \H^1(k_v,{}_gV)^\vee$ via the cup product, the evaluation map
    \[
        \ev \colon {}_gV \otimes {}_gV^\vee \longrightarrow \Fbf_q,
    \]
    and the local invariant. Thus $y_v \in \H^1(k_v,{}_gM^\vee)$ determines a linear map $y_v \colon \H^1(k_v,{}_gV) \to \Fbf_q^d$ via
    \[
        \H^1(k_v,{}_gM^\vee) \simeq \Fbf_q^d \otimes_{\Fbf_q} \H^1(k_v,{}_gV)^\vee \simeq \Hom_{\Fbf_q} \bigl(\H^1(k_v,{}_gV), \Fbf_q^d\bigr).
    \] 
    It follows from the definition of $\psi$ that, under the local invariant isomorphism, $\psi_\ast(x_v\cup y_v)$ corresponds to the composite
    \[
        \Fbf_q^d \overset{x_v}{\longrightarrow}
        \H^1(k_v,{}_gV) \overset{y_v}{\longrightarrow} \Fbf_q^d.
        \]
    This follows from a direct linear algebra calculation. Consequently,
    \[
        \rank(A_v) \leq \dim_{\Fbf_q}\H^1(k_v,{}_gV).
    \]
    We now study the dimension $r_v \coloneq \dim_{\Fbf_q }\H^1(k_v,{}_gV)$. For each finite place $v$ of $k$, write $\Ocal_v$ for the ring of integers of $k_v$. Tate's theorem \cite[\S{10.3}]{Harari2017} and Tate duality \cite[Théorème~10.9]{Harari2017} applied to ${}_gV$ give the ``multiplicative Euler characteristic''
    \[
        \frac{\#\H^0(k_v,{}_gV)\cdot\#\H^2(k_v,{}_gV)}{\#\H^1(k_v,{}_gV)}=\frac{1}{[\Ocal_v : q^n\Ocal_v]}
    \]
    and $\#\H^2(k_v,{}_gV) = \#\H^0(k_v,({}_gV)^D)$, where the superscript $D$ denotes the Cartier dual. Thus we have $r_v = \dim_{\Fbf_q}(\H^1(k_v, {}_gV))\leq n \cdot(2+[k_v:\Qbf_p]) < d$ for places $v \mid p$, and $r_v \leq 2n < d-1$ for all other finite places $v$. At a real place $v$ we have $r_v \leq n < d-1$ because a $1$-cocycle of $\Gal(\Cbf/\Rbf)$ is determined by its value on complex conjugation.

    Thus $\rank(A_v) \leq r_v < d$ at every place $v \mid p$ and $\rank(A_v) < d-1$ at every other finite or real place $v$. If $\lambda_v \neq 0$, equations~\eqref{eq:A_v} and~\eqref{eq:A_vOtherPlaces} would give $\rank(A_v) \geq d-1$ for $v \in \{v_1,\ldots,v_d\}$ and $\rank(A_v) = d$ for $v \notin \{v_1,\ldots,v_d\}$, a contradiction (we recall that $v_1,\ldots,v_d \nmid p$ by construction). Thus $\lambda_v = 0$ for all such $v$, so $A_{v_i} = E_{i,i}$ and $A_v=0$ at every other finite or real place. Applying Lemma~\ref{lemm:GlobalReciprocity} to the global class $\psi_\ast(x\cup y) \in \H^2(k,\Mat_d(\Fbf_q))$, we obtain
    \[
    0 = \sum_{v\ \mathrm{finite\ or\ real}} A_v
    = \sum_{i=1}^d E_{i,i}=I_d,
    \]
    a contradiction.
\end{proof}

To summarize, if $G$ is perfect, $V$ is irreducible and nontrivial, and either $Z(G)=1$ or the action of $G$ on $V$ is faithful, then $Q$ is perfect by Lemma~\ref{lemm:GroupQ}~\ref{lemm:GroupQ:2}. For $d$ as in Proposition~\ref{prop:NoRationalPoint}, the resulting spaces violate the Hasse principle for rational points but satisfy the Hasse principle for zero-cycles of degree $1$ by
Corollary~\ref{coro:HassePerfectStabilizer}.

\begin{lemm}\label{lemm:GlobalReciprocity}
    Let $p$ be a prime number and let $k$ be a number field containing $\mu_p$. Fix an identification $\mu_p \simeq \Fbf_p$, and let $A$ be a finite-dimensional $\Fbf_p$-vector space with trivial Galois action. For each finite or real place $v$ of $k$, we have an isomorphism functorial in $A$
    \[
        \Inv_v^A \colon \H^2(k_v,A) \overset{\simeq}{\longrightarrow} A
    \]
    induced by the local invariant map, using $\frac{1}{p}\Zbf/\Zbf \simeq \Fbf_p$, $a/p \mapsto a \bmod p$, at every place. Then the sequence
    \[
        0\longrightarrow \H^2(k,A) \xrightarrow{\;(\Inv_v^A\circ\operatorname{res}_v)_v\;} \bigoplus_{\substack{v\text{ finite} \\ \text{or real}}} A \xrightarrow{\;\sum_v\;} A \longrightarrow 0
    \]
    is exact.
\end{lemm}
\begin{proof}
    We have $\H^2(k_v,A) = 0$ at complex places, and real places can occur only when $p=2$. Exactness at the last two nonzero terms follows from the nine-term Poitou--Tate sequence, see e.g. \cite[Théorème~17.13]{Harari2017}. The injectivity of the first map follows from Poitou--Tate duality and Chebotarev's density theorem, since the Cartier dual of $A$ is constant.
\end{proof}

\begin{exam} \label{exam:MinimalExampleA5}
    We consider the ``minimal'' example where $p = 2$, $k = \Qbf$, $G = A_5 \cong \SL_2(\Fbf_4)$, $V = \Fbf_4^2$ with the natural action of $G$, and $d=7$. Construct a homogeneous space $X = X_\beta$ as above by choosing the primes $3$, $5$, $7$, $11$, $13$, $17$, and $19$. Its geometric stabilizer is $Q$, where
    \[
        Q = F \rtimes \SL_2(\mathbf F_4), \qquad M = (\Fbf_4^2)^7, \qquad Z = \Mat_7(\Fbf_4)/\Fbf_4 \cdot I_7.
    \]
    Here $F$ is the central extension of $M \oplus M^\vee$ by $Z$ defined by $\phi$ in the construction. Then $|Q| = 15 \cdot 2^{154}$. Thus $X$ has points over every completion of $\Qbf$ and a zero-cycle of degree $1$ over $\Qbf$, even though $X(\Qbf) = \varnothing$.
\end{exam}


\subsection{Explicit zero-cycles of degree one in the absence of rational points} \label{subsec:ExampleNilpotentStabilizers}

In the counterexamples to the Hasse principle on homogeneous spaces with perfect geometric stabilizers constructed in Subsection~\ref{subsec:ExamplePerfectStabilizers}, Corollary~\ref{coro:HassePerfectStabilizer} shows that zero-cycles of degree one do exist there. This statement is not effective. In this subsection we show how to explicitly construct finite extensions of the base field $k$ with mutually coprime degrees, over which the constructed homogeneous spaces have rational points. 

To this end, we give a converse to Proposition~\ref{prop:NoRationalPoint}. First we need the following lemma.

For cup products with coefficients in $\Fbf_p^d$, we recall the identification
\[
    \Fbf_p^d\otimes_{\Fbf_p}\Fbf_p^d\simeq\Mat_d(\Fbf_p),\qquad u\otimes v\longmapsto(v_i u_j)_{1\le i,j\le d}.
\]

\begin{lemm} \label{lemm:LocalDecompositionIntoCupProduct}
    Let $A \in \H^2(k, \Mat_{d}(\Fbf_p))$ be a class for which
    \[
        \rank(\Inv_v^{\Mat_{d}(\Fbf_p)}(A)) \le r_v \coloneq \dim_{\Fbf_p} \H^1(k_v, \Fbf_p)
    \]
    for every finite or real place $v$. Then there exist $x,y \in \H^1(k,\Fbf_p^d) = \H^1(k,\Fbf_p)^d$ such that
    \[
        x \cup y = A.
    \]
\end{lemm}
\begin{proof}
    At complex places, take $x_v=y_v=0$. At every finite or real place $v$, the cup-product pairing
    \[
        \H^1(k_v, \Fbf_p) \times \H^1(k_v, \Fbf_p) \longrightarrow \H^2(k_v,\Fbf_p)
    \]
    is a perfect duality of $r_v$-dimensional $\Fbf_p$-vector spaces, by Tate's local duality. Thanks to the assumption on the rank of the local invariants of $A$, we may write
    \[
        A|_{k_v} = x_v \cup y_v,
    \]
    where $x_v,y_v \in \H^1(k_v,\Fbf_p^d) = \H^1(k_v,\Fbf_p)^d$. Now let $S \subseteq \Omega_k$ be a finite subset containing the archimedean places, the places lying above $p$, and the places $v$ for which $A|_{k_v} \neq 0$. Recall that
    \[
        \H^1(k, \Fbf_p^d) \cong (k^\times / k^{\times p})^d \qquad \text{and} \qquad \H^1(k_v, \Fbf_p^d) \cong (k_v^\times / k_v^{\times p})^d
    \]
    by Kummer's theory. The Artin--Whaples approximation lemma therefore gives an element $x \in \H^1(k,\Fbf_p^d)$ for which $x|_{k_v} = x_v$ for all $v \in S$, whence $x|_{k_v} \cup y_v = A|_{k_v}$. For $v \notin S$ we have $x|_{k_v} \cup 0 = 0 = A|_{k_v}$. Then the ``arithmetic lemma'' in \cite[Proposition~5.2]{Linh2024BK} provides an element $y \in \H^1(k,\Fbf_p^d)$ for which $x \cup y = A$.
\end{proof}

Retain all notation from Subsection~\ref{subsec:ExamplePerfectStabilizers}, but {\bf drop the assumption} 
\[
    d > \max\limits_{v \mid p} n([k_v : \Qbf_p] + 2).
\]

\begin{prop} \label{prop:RationalPointCriterion}
    For a homogeneous space $X = X_\beta$ constructed from a class $\beta \in \H^2(k, Z)$ as in Subsection~\ref{subsec:ExamplePerfectStabilizers}, the following are equivalent.
    \begin{enumerate}
        \item \label{prop:RationalPointCriterion:1} There is a place $v_0 \mid p$ with $d \le n([k_{v_0} : \Qbf_p] + 2)$.
    
        \item \label{prop:RationalPointCriterion:2} There exist $x \in \H^1(k,M)$ and $y \in \H^1(k,M^\vee)$ for which $\beta = \phi_\ast(x \cup y)$.
        
        \item \label{prop:RationalPointCriterion:3} $X(k) \neq \varnothing$.
    \end{enumerate}
    Moreover, let $K/k$ be a finite extension. If there exists a place $w_0\mid p$ of $K$ such that $d\le n([K_{w_0}:\Qbf_p]+2)$, then $X(K)\neq\varnothing$.
\end{prop}
\begin{proof}
    Assume~\ref{prop:RationalPointCriterion:1}. Since $v_0 \notin \{v_1,\ldots,v_d\}$, we may define $A \in \H^2(k, \Mat_{d}(\Fbf_p))$ by
    \[
        \begin{cases}
            \Inv_{v_i}^{\Mat_{d}(\Fbf_p)}(A) = E_{i,i} & \text{for } i = 1,\ldots,d,\\
            \Inv_{v_0}^{\Mat_{d}(\Fbf_p)}(A) = - I_d,\\
            A|_{k_v} = 0 & \text{for } v \notin \{v_0,v_1,\ldots,v_d\}.
        \end{cases}
    \]
    Indeed, this is possible thanks to Lemma~\ref{lemm:GlobalReciprocity}. Clearly $\beta = \pi_\ast(A)$. Next, write $\{1,\ldots,d\}$ as a disjoint union of $n$ (possibly empty) subsets $J_1,\ldots,J_n$
    of size $|J_j| \le [k_{v_0} : \Qbf_p] + 2$. Again Lemma~\ref{lemm:GlobalReciprocity} yields an element $B_j \in \H^2(k, \Mat_{d}(\Fbf_p))$ with 
     \[
    \begin{cases}
        \Inv_{v_i}^{\Mat_{d}(\Fbf_p)}(B_j) = E_{i,i} & \text{for } i \in J_j,\\
        \Inv_{v_0}^{\Mat_{d}(\Fbf_p)}(B_j) = - \sum_{i \in J_j} E_{i,i},\\
        B_j|_{k_v} = 0 & \text{for } v \notin \{v_0\} \sqcup \{v_i \mid i \in J_j\}.
    \end{cases}
    \]
    Then $\sum_{j = 1}^n B_j = A$. Furthermore $\rank(\Inv_v^{\Mat_d(\Fbf_p)}(B_j)) \le \dim_{\Fbf_p} \H^1(k_v, \Fbf_p)$ for every finite or real place $v$ (at $v = v_0$ this follows from \cite[Chapter~VII, Proposition~6.8]{Milne2020CFT}, and at the other places this is trivial). By Lemma~\ref{lemm:LocalDecompositionIntoCupProduct} there exist $x_j, y_j \in \H^1(k, \Fbf_p^d) = \H^1(k, \Fbf_p)^d$ with 
    \[
        x_j \cup y_j = B_j.
    \]
    In what follows we abusively identify $x_j$ and $y_j$ with their respective images in $\H^1(k,\Fbf_q)^d$. Choose an $\Fbf_q$-basis $e_1,\ldots,e_n$ for $V$ and denote by $e_1^\vee,\ldots,e_n^\vee$ its dual basis for $V^\vee$. Set
    \[
        x \coloneq \sum_{j=1}^n x_j e_j \in \H^1(k, M) \qquad \text{and} \qquad y \coloneq \sum_{j=1}^n y_j e_j^\vee \in \H^1(k, M^\vee).
    \]
    Then $\psi_\ast(x \cup y) = \sum_{j=1}^n (x_j \cup y_j) = \sum_{j=1}^n B_j = A$, thus 
    \[
        \phi_\ast(x \cup y) = \pi_\ast(A) = \beta, 
    \]
    which proves~\ref{prop:RationalPointCriterion:2}. Now~\ref{prop:RationalPointCriterion:2}$\Rightarrow$\ref{prop:RationalPointCriterion:3} follows from lemmas~\ref{lemm:ConnectingMap} and~\ref{lemm:HomoSpacesCoho}, while~\ref{prop:RationalPointCriterion:3}$\Rightarrow$\ref{prop:RationalPointCriterion:1} is simply the contraposition of Proposition~\ref{prop:NoRationalPoint}.

    For the final assertion, put $m=[K:k]$ and repeat the construction over $K$, prescribing $[K_w:k_{v_i}]E_{i,i}$ at each place $w\mid v_i$ and $-mI_d$ at $w_0$, and making the corresponding replacements for the classes $B_j$. These prescriptions satisfy reciprocity since $\sum_{w\mid v_i}[K_w:k_{v_i}]=m$, and their images under $\pi$ give $\beta|_K$. Choose the blocks with $|J_j|\le [K_{w_0}:\Qbf_p]+2$. The same rank bounds and Lemma~\ref{lemm:LocalDecompositionIntoCupProduct}, applied over $K$, give $\beta|_K=\phi_\ast(x\cup y)$, hence $X(K)\neq\varnothing$.
\end{proof}

\begin{exam} \label{exam:MinimalExampleA5:ExplicitZeroCycle}
    Recall the homogeneous space $X = X_\beta$ from Example~\ref{exam:MinimalExampleA5} with the choices $p = 2$, $k = \Qbf$, $G = A_5$, $n = 2$, $d=7$, and the seven primes $3$, $5$, $7$, $11$, $13$, $17$, and $19$. Its geometric stabilizer is $Q$, where
    \[
        Q = F \rtimes \SL_2(\mathbf F_4), \qquad M = (\Fbf_4^2)^7, \qquad Z = \Mat_7(\Fbf_4)/\Fbf_4 \cdot I_7.
    \]
    Let $K_2 \coloneq \Qbf(\sqrt{2})$. It is an Eisenstein extension of $\Qbf$ at $p = 2$, so the prime $2$ is totally ramified in $K_2/\Qbf$. Hence the only place $w_2 \in \Omega_{K_2}$ lying above $2$ satisfies $2([(K_2)_{w_2} : \Qbf_2] + 2) = 8 > 7$. On the other hand, consider $K_3 \coloneq \Qbf(\theta)$, where $\theta^3 + \theta + 1 = 0$. Since $T^3+T+1$ is irreducible modulo $2$, the prime $2$ remains inert in $K_3/\Qbf$. Thus the only place $w_3 \in \Omega_{K_3}$ lying above $2$ satisfies $2([(K_3)_{w_3} : \Qbf_2] + 2) = 10 > 7$ as well. By Proposition~\ref{prop:RationalPointCriterion} we have $X(K_2) \neq \varnothing$ and $X(K_3) \neq \varnothing$. Thus, $[P_3]-[P_2]$ is a zero-cycle of degree $1$ on $X$, where $P_i$ denotes the closed point defined by any point of $X(K_i)$ for $i=2,3$.
\end{exam}

\begin{exam}
    We specialize the construction of Subsection~\ref{subsec:ExamplePerfectStabilizers} to $G = 1$, $q = p$, $V = \Fbf_p$, and $d \geq [k:\Qbf]+3$; thus $Q = F$ (which is no longer perfect, but of nilpotency class two). Again choose $d$ distinct finite places $v_1,\ldots,v_d$ of $k$ not lying over $p$, and define $\beta$ by the same local invariant prescription as in Subsection~\ref{subsec:ExamplePerfectStabilizers}, with zero invariants at all other places. Such a class exists by Lemma~\ref{lemm:GlobalReciprocity}. Set $X = X_\beta$ as in Lemma~\ref{lemm:HomoSpacesCoho}. Then we obtain the construction of \cite[\S2]{Linh2026Insufficiency}. While $F$ is not perfect, the arguments in \cite[\S\S{2 and 3}]{Linh2026Insufficiency} still show that $X$ violates the Hasse principle for rational points and satisfies
    \[
        \Br_{\nr,1}(X) = \Br_0(X) \qquad \text{and} \qquad \Br_{\nr}(X)/\Br_0(X) \cong \Zbf/p\Zbf.
    \]
    Moreover, its failure of the Hasse principle can be explained by a transcendental Brauer--Manin obstruction. Nevertheless, local points define local zero-cycles of degree one, and since $\Br_{\nr,1}(X)=\Br_0(X)$, Theorem~\ref{thm:Main} ensures the existence of a zero-cycle of degree one on $X$.
    \medskip
    
    \noindent The ``minimal'' example obtained in this fashion is when $p = 2$, $k = \Qbf$, and $d = 4$, with the four primes  $v_1=3$, $v_2=5$, $v_3=7$, and $v_4=11$.
    The geometric stabilizer is
    \[
        F= \left(\Mat_4(\Fbf_2)/\Fbf_2 \cdot I_4 \right) \times_{\phi} \left(\Fbf_2^4 \oplus (\Fbf_2^4)^\vee \right).
    \]
    Here $\times_\phi$ denotes the product endowed with the group law defined in Subsection~\ref{subsec:ExamplePerfectStabilizers}. It has order $|F| = p^{d^2 + 2d-1} = 2^{23} = 8388608$. Again in this example we have $X(K_2) \neq \varnothing$ and $X(K_3) \neq \varnothing$ thanks to Proposition~\ref{prop:RationalPointCriterion}, where $K_2 = \Qbf(\sqrt{2})$ and $K_3 = \Qbf(\theta)$, $\theta^3 + \theta + 1 = 0$. Thus, $[P_3]-[P_2]$ is a zero-cycle of degree $1$ on $X$, where $P_i$ denotes the closed point defined by any point of $X(K_i)$ for $i=2,3$.
\end{exam}

\appendix
\refstepcounter{section}
\section*{Appendix: Compatibility with central extensions}
\label{app:BrauerCompatibility}

\phantomsection
\addcontentsline{toc}{section}{Appendix: Compatibility with central extensions}

\begin{lemm}
\label{lemm:BrauerCentralExtensionCompatibility}
    With the notation of the central-extension construction in
Subsection~\ref{subsec:BrauerGroup}, assume that $\Pic(G_{\bar{K}})=0$. Then the diagram
    \[
    \begin{tikzcd}[column sep=large, row sep=large]
        \H^2(F,\bar{K}^\times) \arrow[r, "\sim"] \arrow[rd, "j"'] & \Ext^{c}_{\bar{K}\gp} (F_{\bar{K}},\Gbf_m) \arrow[d, "\rho_{\tors(G_{\bar{K}})}"] \\
        & \Br(Y_{\bar{K}})
    \end{tikzcd}
    \]
    commutes, where $j$ is the homomorphism of
Proposition~\ref{prop:BrNrSubgroupOfH2}.
\end{lemm}

\begin{proof}
    Since the three morphisms come from three different articles, we will draw a large diagram with all of the groups concerned in these three definitions.

    Let $\alpha\in \H^2(F,\bar{K}^\times)$, and let
    \[
        e_\alpha \colon \qquad 1 \longrightarrow \Gbf_m \longrightarrow E_\alpha \longrightarrow F_{\bar{K}} \longrightarrow 1  
    \]
    be the corresponding central extension. We prove that
    \[
        j(\alpha) = \rho_{\tors(G_{\bar{K}})}([e_\alpha]) \quad \text{in } \Br(Y_{\bar{K}}).
    \] 
    Since $Y_{\bar{K}}$ is smooth and integral, restriction to its generic point $\eta$ induces an injection
    \[
       \eta^\ast \colon \Br(Y_{\bar{K}}) \longhookrightarrow \Br(\bar{K}(Y)).
    \]
    It therefore suffices to compare the restrictions of the two classes to the generic point.

    Fix an embedding of $\bar{K}(G)$ into an algebraic closure $\overline{K(Y)}$, and write
    \[
        \Gamma_\eta \coloneq \Gal\bigl(\overline{K(Y)} / \bar{K}(Y)\bigr).
    \]
    The generic fibre $\Spec(\bar{K}(G)) \to \Spec(\bar{K}(Y))$ is an $F_{\bar{K}(Y)}$-torsor, represented by a continuous homomorphism $\rho_\eta \colon \Gamma_\eta \twoheadrightarrow F$. We use the corresponding identification $\Gal(\bar{K}(G)/\bar{K}(Y)) \simeq F$. The constructions under consideration fit into the following diagram, whose commutativity we establish below:
    \[
    \begin{tikzcd}[column sep=large, row sep=huge,
    cells={nodes={text height=2.5ex, text depth=0.75ex}}]
        \H^2(F, \bar{K}^\times) \arrow[r, "{\substack{\text{change of} \\  \text{coefficients}}}"]
        \arrow[d, equal] & \H^2(F, \bar{K}[G]^\times)
        \arrow[r] \arrow[d, "{\substack{\text{change of}\\ \text{coefficients}}}"] & \Br(Y_{\bar{K}}) \arrow[d, "\eta^\ast"] \\
        \H^2(F, \bar{K}^\times) \arrow[r, "{\substack{\text{change of}\\ \text{coefficients}}}"] \arrow[d, "\rho_\eta^\ast"'] & \H^2(F, \bar{K}(G)^\times) \arrow[r, "\inf"] \arrow[d, "\inf"] & \Br(\bar{K}(Y)) \arrow[d, equal] \\
        \H^2(\Gamma_\eta, \bar{K}^\times) \arrow[r, "{\substack{\text{change of}\\ \text{coefficients}}}"] & \H^2\!\left(\Gamma_\eta, \overline{K(Y)}^\times\right) \arrow[r, equal] & \Br(\bar{K}(Y)) \\
        \H^2(F,\bar{K}^\times) \arrow[r] \arrow[u, "\rho_\eta^\ast"] & \Ext^{c}_{\bar{K}(Y)\gp} (F_{\bar{K}(Y)}, \Gbf_m) \arrow[r, "\rho_{\tors,\eta}"] \arrow[u, "\rho_{\tors,\eta}"'] & \Br(\bar{K}(Y)) \arrow[u, equal]
        \\
        \H^2(F, \bar{K}^\times) \arrow[r, "\sim"] \arrow[u, equal] &
        \Ext^{c}_{\bar{K}\gp} (F_{\bar{K}}, \Gbf_m) \arrow[r, "\rho_{\tors(G_{\bar{K}})}"'] \arrow[u, "\text{base change}"'] & \Br(Y_{\bar{K}}) \arrow[u, "\eta^\ast"']
    \end{tikzcd}
    \]
    The maps in the first row are those constructed in the proof of Proposition~\ref{prop:BrNrSubgroupOfH2}; their composite is $j$. The arrows labeled $\inf$ denote inflation to $\Gamma_\eta$ followed by the map induced by the inclusion $\bar{K}(G)^\times \hookrightarrow \overline{K(Y)}^\times$. The first arrow in the fourth row sends $\alpha$ to $[(e_\alpha)_{\bar{K}(Y)}]$. Here $\rho_{\tors,\eta}$ denotes the homomorphism of \cite[\S~2, pp.~313--314]{CTX09} associated with the generic torsor $\Spec(\bar{K}(G)) \to \Spec(\bar{K}(Y))$, taking $[e]$ to $\partial_e([\rho_\eta])$.

    By functoriality of the Hochschild--Serre spectral sequence, the first two rows are compatible with restriction to the generic point. The left-hand square between the second and third rows commutes by compatibility of inflation with change of coefficients, which is immediate at the level of cocycles. The right-hand square between the last two rows commutes by functoriality in the base scheme of the pairing in \cite[\S~2, p.~314]{CTX09}, applied to $\eta$.

    It remains to compare the two paths in the left-hand square between the third and fourth rows. Choose a normalized cocycle $c$ representing $\alpha\in\H^2(F,\bar{K}^\times)$. Along the first path, pullback by $\rho_\eta$ followed by the inclusion $\bar{K}^\times \hookrightarrow \overline{K(Y)}^\times$ gives the class of the cocycle
    \[
        (\sigma, \tau) \longmapsto c \bigl(\rho_\eta(\sigma), \rho_\eta(\tau)\bigr) \in \bar{K}^\times \subseteq\overline{K(Y)}^\times.
    \]
    Along the other path, we represent $e_\alpha$ by the extension whose middle term has underlying scheme $\Gbf_m \times F_{\bar{K}}$ and multiplication given by
    \[
        (u,s)(v,t) \coloneq (uv\,c(s,t),st).
    \]
    Base change to $\bar{K}(Y)$ gives the central extension $(e_\alpha)_{\bar{K}(Y)}$ of $F_{\bar{K}(Y)}$ by $\Gbf_m$, with the same multiplication law. Taking $\overline{K(Y)}$-points gives a central extension of $\Gamma_\eta$-groups
    \[
        1 \longrightarrow \overline{K(Y)}^\times  \longrightarrow E_\alpha \Bigl(\overline{K(Y)}\Bigr) \longrightarrow F \longrightarrow 1.
    \]
    Here the action of the Galois group is given by $\sigma(u,s) \coloneq (\sigma(u),s)$. To compute the boundary, choose the lift $t_s \coloneq (1,s) \in E_\alpha(\overline{K(Y)})$ for each $s \in F$. Then $\sigma(t_s) = t_s$ for every $\sigma \in \Gamma_\eta$. Since $\rho_\eta$ is a homomorphism, the boundary class $\partial_{(e_\alpha)_{\bar{K}(Y)}}([\rho_\eta])$ is represented by the $2$-cocycle
    \begin{align*}
        \Gamma_\eta \times \Gamma_\eta \longrightarrow \overline{K(Y)}^\times, \qquad (\sigma,\tau) & \longmapsto t_{\rho_\eta(\sigma)}\,\sigma\bigl(t_{\rho_\eta(\tau)}\bigr)
        \, t_{\rho_\eta(\sigma\tau)}^{-1}\\
        & = t_{\rho_\eta(\sigma)}\,t_{\rho_\eta(\tau)} \, t_{\rho_\eta(\sigma)\rho_\eta(\tau)}^{-1}\\
        & = c \bigl(\rho_\eta(\sigma), \rho_\eta(\tau)\bigr).
    \end{align*}
    Thus the two paths give the same cocycle in $\overline{K(Y)}^\times$, and the square commutes. The two classes have the same restriction to the generic point, and hence are equal by injectivity of $\eta^\ast$.
\end{proof}

{\em AI Declaration.} In certain lemmas involving commutative diagrams, we used ChatGPT to produce portions of the proofs of commutativity. The arguments provided by ChatGPT have been understood, and the relevant proofs were substantially rewritten for clarity. We also used ChatGPT to proofread the manuscript and revise the exposition. The entire article has been checked line by line by the authors.
\medskip

{\em Acknowledgments.} Nguy\~{\^{e}}n M\d{a}nh Linh acknowledges support from the Fondation Sciences Math{\'e}matiques de Paris during his postdoctoral appointment. He is grateful to his Ph.~D. supervisor, David Harari, for raising this question during his defense.
Tianyang Wang sincerely thanks his Ph.~D. supervisor, Cyril Demarche, for valuable discussions throughout this work and continued support during his doctoral studies, and gratefully acknowledges the financial support provided by Sorbonne Université through a doctoral contract. We both appreciate the Institut de Math{\'e}matiques de Jussieu - Paris Rive Gauche for excellent working conditions.

\bibliographystyle{alpha}
\bibliography{ref}
\end{document}